\documentclass[a4paper,reqno]{amsart}

\providecommand\nonumsection{\section*}

\usepackage{miscmath} 
\numberwithin{equation}{section} 
\usepackage[ps,dvips,xdvi,arc,all,knot,frame,curve,poly]{xy}
\usepackage{rg} 

\usepackage{rotating}

\usepackage{graphicx}
\newcommand\ignoregroup[1]{}
\newcommand\usegraphics[2][0pt]{%
  \IfFileExists{#2}%
    {\raisebox{#1}{\includegraphics{#2}}\ignoregroup}%
    {}%
}

\usepackage{prettyref}
\newrefformat{appendix}{Appendix~\ref{#1}}
\newrefformat{prop}{Proposition~\ref{#1}}
\newrefformat{lemma}{Lemma~\ref{#1}}
\newrefformat{xmp}{Example~\ref{#1}}
\newrefformat{dfn}{Definition~\ref{#1}}
\newrefformat{cor}{Corollary~\ref{#1}}
\newrefformat{rem}{Remark~\ref{#1}}

\usepackage[newitem,newenum]{paralist}
\defaultenum{(i)}{(a)}{(1)}{(A)}

\newcommand{\catA}{\category{A}}     
\newcommand{\catN}{\category{N}}     
\newcommand{\catVect}{\category{Vect}}
\newcommand{\EndOp}[1][]{\operad{E}_{#1}} 
\newcommand{\catSet}{\category{Set}} 
\newcommand{\Oo}[1][]{\operad{O}_{#1}}  
\newcommand{\Pp}[1][]{\operad{P}_{#1}}  
\newcommand{\Dia}[1][]{\operad{D}_{#1}}  
\newcommand{\RT}[1][]{\operad{T}_{#1}}  
\newcommand{\RG}[1][]{\operad{R}_{#1}}  
\newcommand{\CG}[1][]{\operad{G}_{#1}}  
\newcommand{\freemc}[1]{#1^{\otimes}}        
\newcommand{\freemsc}[1]{#1^{\bigstar}}     
\newcommand{\rdual}[1]{{#1}\spcheck}
\newcommand{\ldual}[1]{{}\spcheck{#1}}
\newcommand{\opp}{\sptext{opp}}

\DeclareMathOperator{\Src}{Src}
\DeclareMathOperator{\Tgt}{Tgt}
\DeclareMathOperator{\In}{In}
\DeclareMathOperator{\Out}{Out}
\DeclareMathOperator{\Leg}{Leg}

\newcommand{\correlator}[1]{\langle #1 \rangle}

\newcommand{\M}{\moduli} 
\newcommand{\Mbar}{\modulibar}

\newcommand{\Hermitian}[1]{{\mathcal H}}

\newcommand{\Vertices}[1]{#1^{(0)}}
\newcommand{\Edges}[1]{#1^{(1)}}
\newcommand{\Holes}[1]{#1^{(2)}}

\newcommand{\avg}[2][V]{\langle #2\ \!\rangle_{#1}^{}} 
\newcommand{\cplx}[1]{#1^{}_\setC} 

\begin{document}

\title{Feynman Diagrams via Graphical Calculus}
\date{2001-06-23}
\author{Domenico Fiorenza}
\address[D. Fiorenza]{%
    Dipartimento di Matematica ``Leonida Tonelli'' \\
    Universit{\`a} di Pisa, \\
    via F. Buonarroti, 2 \\
    56127 Pisa \\
    Italy
    }
\email{fiorenza@mail.dm.unipi.it}

\author{Riccardo Murri}
\address[R. Murri]{%
  Scuola Normale Superiore \\
  p.za dei Cavalieri, 7 \\
  56127 Pisa \\
  Italy
  }
\email{murri@cibsmail.sns.it}

\renewcommand{\subjclassname}{\textup{2000} Mathematics Subject
     Classification}

\maketitle

\setcounter{tocdepth}{1} 
\tableofcontents

\section{Introduction}
\label{sec:introduction}

It has been known for a while, in Hopf algebraists' folklore, that
there is a very close connection between the graphical formalism for
ribbon categories and Feynman diagrams.
Although this correspondence is frequently implied, it seems to have been
first explicitly described in the recent \cite{oeckl;braided-qft}. 
Yet, we know of no systematic exposition in existing literature; the
aim of this paper is to provide such an account. In particular, in
deriving Feynman diagrams expansion of Gaussian integrals as an
application of the
graphical formalism for symmetric monoidal categories, we discuss in
detail how different kinds of interactions give rise to different families 
of graphs and show how symmetric and cyclic interactions lead to
``ordinary'' and ``ribbon'' graphs respectively.

Feynman diagrams are usually introduced as a kind of combinatorial
bookkeeping device in asymptotic expansion of Gaussian integrals
(cf. \cite{bessis-itzykson-zuber;graphical-enumeration}). Indeed,
given any integral 
\begin{equation}\label{eq:GI}
  \int_{\mathcal{H}} f(X) \E^{-S(X)} \ud \mu(X),
\end{equation}
(where $f$ and $S$ are polynomial functions, $\ud
\mu$ the Gaussian measure over a real Hilbert space $\mathcal{H}$) its
asymptotic expansion can be written in terms of ``correlator
functions''
\begin{equation}
  \label{eq:SC}
  \correlator{X_1 \cdots X_k} := \int_{\mathcal{H}} X_1 \cdots X_k \ud \mu(X),
\end{equation}
where $X_1, \ldots, X_k$ are coordinates of $X$ with respect to a chosen
basis of $\mathcal{H}$. The data identifying correlators can be put in
a one-to-one correspondence with some combinatorial data describing a 
graph; vice-versa, a correlator may be reconstructed from a graph by
some simple ``Feynman rules''. Therefore, the asymptotic expansion of
\prettyref{eq:GI} can be written as a sum over graphs.

In this paper we take the other way round: we associate an analytic
expression to a graph by means of graphical calculus, then we show that
the summation of all these expressions, for some chosen class of
graphs, gives the asymptotic expansion of an integral of the type
\prettyref{eq:GI}. A relation between the integrand $f$ in
\prettyref{eq:GI} and the class of graphs being summed upon is
derived, and found to coincide with usual ``Feynman rules''. This
point of view has an advantage: one can easily reconstruct a ``path
integral'' formulation from a set of given Feynman rules.

\prettyref{sec:gc} presents a sketchy account of graphical calculus,
in its ``ribbon graphs with coupons'' flavor, following mainly
\cite{reshetikhin-turaev;ribbon-graphs}. We assume the reader is
already acquainted with the material presented there, and do not
provide any proof nor motivation for this theory. A very readable
introduction to this sort of graphical calculus may be found in
\cite{bakalov-kirillov}; for a comprehensive treatment and further
references consult \cite{joyal-street;1991, kassel;quantum-groups,
turaev;invariants-of-knots}.

We specialize graphical calculus to the vector spaces category
(equipped with the trivial braiding); this would be a very
uninteresting choice in the context of knot theory --- where graphical
calculus was originally developed --- since it gives rise to trivial
link invariants. However, it is the right way to go here, because we
want the analytic expression associated with a diagram to depend only
on the topology of the diagram, and not on its immersion in the plane.

In \prettyref{sec:feynman-diagrams} the relevant theorems relating
graphical calculus and Feynman diagrams expansion of Gaussian
integrals are stated and proved. It is ``folklore'' material, and we
know of no other written reference for it.

\prettyref{sec:matrix-models} works out in detail an example of
physical interest, namely, the Kontsevich model for 2D quantum gravity
(of which the 't~Hooft standard matrix model is a particular case),
showing how the new notation can be applied to known cases.

Gaussian integrals and Feynman diagrams have been generalized by Robert
Oeckl to the wider context of (not necessarily symmetric) braided monoidal
categories. In particular, using graphical calculus techniques, he proves
that any braided Gaussian integral admits an expansion in braided Feynman
diagrams (see \cite{oeckl;braided-qft} for details). 

\begin{notations}\label{notations:notations}
If $\catA$ is a category, we write $X \in \catA$ to state that $X$ is
an object of $\catA$. The map notation $f: X \to Y$ will be
occasionally used to denote a morphism $f \in \catA(X,Y)$.

The symbol $\Perm{k}$ stands for the permutation group on $k$ letters.

There are several classes of graphs appearing in the text: we have
reserved the term ``ordinary'' graph for purely $1$-dimensional
CW-complexes; all other graphs (ribbon, RT, modular) differ by some
additional structure on the vertices --- precise definitions follow in
the body of the text.

Unfortunately, there seems to be no agreement among authors about the
naming of objects involved in graphical calculus; our own choice, to
the readers' bewilderment, is a mixture of many naming styles found in
the literature, and is not entirely consistent with any of our
sources.
\end{notations}


\everyxy={/r24pt/:} 

\section{Feynman diagrams via graphical calculus}
\label{sec:gc}

In this section we recall some basic facts of graphical calculus, as
introduced by Reshetikhin-Turaev in
\cite{reshetikhin-turaev;ribbon-graphs} and Joyal-Street in
\cite{joyal-street;1991}; in particular, we state the main result of
this theory in the simpler case of ribbon graphs (in the sense of
\cite{kontsevich;intersection-theory;1992}).  We make a fundamental
simplification in our exposition of this theory, namely, we drop the
requirement that graphs edges are equipped with a framing: indeed,
since our ground category has trivial \emph{balancing} (cf.
\cite{joyal-street;1991} and \cite{shum;tortile-categories}), we do
not need the extra structure given by twists.  We stick to the usual
``wireframe'' graphs, which somewhat simplifies definitions.

\subsection{Preliminaries on tensor categories and PROPs}
\label{sec:tensor-categories}

The notion of monoidal category is well-known and discussed at length
in the existing literature; we recall a few facts and definitions. A
precise list of axioms may be found in the already cited sources.

A monoidal category \((\catA, \otimes, I, a, l, r)\) is given by
the following data:
\begin{enumerate}
\item a category \(\catA\);
\item a functor \(\otimes: \catA \times \catA \to \catA\);
\item a distinguished object (identity) \(I\in\catA\);
\item a natural transformation  \(a: \bigl( (- \otimes -) \otimes - \bigr)
\to
\bigl( - \otimes (- \otimes -) \bigr)\) of functors \(\catA \times \catA
\times \catA
\to \catA\);
\item a natural transformation \(l: (- \otimes I) \to \Id_{\catA}\);
\item a natural transformation \(r: (I \otimes -) \to \Id_{\catA}\).
\end{enumerate}
These data are required to satisfy some compatibility axioms: roughly
speaking, one would consider \(\otimes\) as a ``multiplication'' of
objects in the category \(\catA\), and \(a\), \(l\), \(r\) are the
appropriate ``categorizations'' of usual conditions expressing
associativity of multiplication and existence of a bilateral
multiplicative identity \(I\). Indeed, these axioms imply that
expressions like \(X_1 \otimes X_2 \otimes \dots \otimes X_k\) are
well-defined (i.e., do not depend on the way we put parentheses in them)
up to a natural isomorphism, and that insertion or removal of \(I\)
can be neglected, again up to a natural isomorphism.

If \(\catA\) is an Abelian category we require \(\otimes\), \(a\),
\(l\), \(r\) to be linear. Abelian monoidal categories are called
\emph{tensor} categories.

A tensor functor \(F\) is a functor such that \(F(A \otimes B) =
F(A) \otimes F(B)\), up to a natural isomorphism.

A braided tensor category is a tensor category equipped with a
family of isomorphisms \(\tau_{XY}: X \otimes Y \to Y \otimes X\),
natural in both \(X\) and \(Y\), satisfying a certain compatibility
diagram (MacLane's hexagon condition). In a braided tensor category
any two expressions \(X_1 \otimes X_2 \otimes \dots \otimes X_r\) and
\(X_{\sigma_1} \otimes X_{\sigma_2} \otimes \dots X_{\sigma_r}\) (with
\(\sigma\in\Perm{k}\)) are isomorphic, but there are (possibly) many
different isomorphisms built from maps \(\tau_{X_iX_j}\).

A symmetric tensor category \(\catA\) is a braided category such
that \(\tau_{YX} \circ \tau_{XY} = \id_{X\otimes Y}\) for all \(X\), \(Y\)
objects of \(\catA\). This implies that an isomorphism between \(X_1
\otimes X_2 \otimes \dots \otimes X_r\) and \(X_{\sigma_1} \otimes
X_{\sigma_2} \otimes \dots X_{\sigma_r}\), built only from maps
\(\tau_{X_i X_j}\), depends only on the permutation \(\sigma\).

A braided tensor category \(\catA\) is said to have right duals if
for any object \(X\in\catA\) there is an object \(\rdual{X}\) and
morphisms \(\ev_X: X \otimes \rdual{X} \to \fk\) and \(\coev_X: \fk \to
\rdual{X} \otimes X\) such that the composition
\begin{equation*}
X \xrightarrow{\id \otimes \coev} X \otimes \rdual{X} \otimes X
\xrightarrow{\ev \otimes \id} X
\end{equation*}
is the identity morphism on $X$.  Definition of left duals is
completely analogous.  A braided tensor category is rigid if it has
both left and right duals and they are canonically isomorphic. For any
object $A$ of a rigid tensor category $\catA$, put
\begin{equation*}
A^r := 
\begin{cases}
A\tp{r}, \quad        &\text{if }r>0, \\
I,                    &\text{if }r=0, \\
(\rdual{A})\tp{(-r)}, &\text{if }r<0.
\end{cases}
\end{equation*}
One may check that the usual relation $A^r \otimes A^s = A^{r+s}$ holds,
up
to a natural isomorphism.

\begin{xmp}
The category of vector spaces, equipped with the usual tensor product
and the obvious $a$, $l$, $r$ is rigid and symmetric.
\end{xmp}

\subsubsection{Free tensor categories}
\label{sec:free-mono-categories}

For any category \(\catA\) we can form a monoidal category
\(\freemc{\catA}\): objects are finite sequences of objects in
\(\catA\), and morphisms are finite sequences of morphisms from
\(\catA\). Tensor product is given by juxtaposition; the identity
object is the empty sequence; associativity and identity natural
transformations are the obvious ones. If \(\catA\) is itself monoidal,
then there is an obvious functor \(\freemc\catA \to \catA\).

\subsubsection{An important symmetric rigid tensor category}
\label{sec:vector-space-cat}
Fix a vector space \(V\) (over a ground field \(\fk\)) and a
non-degenerate symmetric inner product \(b: V\otimes V\to\fk\) on it.
Build a category
\(\langle V \rangle\): objects are powers \(V\tp{k}\) of \(V\), for
$k\in\setN$, and
\(I:=\fk\simeq V\tp{0}\); morphisms are linear maps \(V\tp{r} \to
V\tp{s}\).  This category is symmetric with the usual tensor
product of vector spaces.

The object \(V\) is self-dual, if we define the morphism \(\ev_V\) to be
the inner product \(b: V\otimes V \to \fk\), and \(\coev_V\) to be the
morphism sending \(1\in\fk\) to the Casimir element \(\sum e_i\otimes
e^i\), where \(\{e_i\}\) is a basis of \(V\) and \(\{e^i\}\) is the
dual one. Similarly, one can define morphisms such that any
\(V\tp{r}\) is left and right dual to itself.  As a result, \(\langle
V \rangle\) is rigid.

\subsection{PROPs}
\label{sec:props}
Informally speaking, a PROP is a monoidal category whose
\(\Hom\)-spaces are objects of another monoidal category: if $\catA$
is a category, then, for any two objects $X,Y\in \catA$, $\catA (X,Y)
:= \Hom_{\catA} (X, Y)$ is a \emph{set}, what is more, $\Hom:
\catA\opp \times \catA \to \catSet$ is a functor; a structure of
$\catA_{\text{Hom}}$-PROP $\Pp$ over $\catA_{\text{Ob}}$ is given by a
functor $\Pp: (\catA_{\text{Ob}})\opp \times \catA_{\text{Ob}} \to
\catA_{\text{Hom}}$ whose properties generalize those of the
$\Hom$-functor. A precise definition of PROP is rather cumbersome, so
we consign it to \prettyref{appendix:PROPs}; here we shall give only
some illustrative examples.

If \(X\) and \(Y\) are vector spaces, then \(\Hom(X,Y)\) is a vector
space: if we define $\operad{V}(X, Y) := \Hom(X,Y)$, then $\operad{V}$
is a natural structure of a \(\catVect\)-PROP over the category
\(\catVect\).

Similarly, if \(H_1\) and \(H_2\) are Hilbert spaces,
\(\Hom(H_1,H_2)\) is a Banach space: the category of Hilbert spaces is
in a natural way a Banach spaces-PROP.

Likewise, every tensor category is a $\catVect$-PROP with the trivial
PROP structure given by $\Pp(X, Y) := \Hom(X, Y)$; by abuse of
language we say that $\Pp(X, Y)$ are the $\Hom$-spaces of the PROP.
The PROPs of graphs we are going to introduce are indeed tensor
categories, but we are mainly interested in their actions as PROPs.

\subsubsection{PROP-algebras}
\label{sec:prop-algebras}
Now, build a category $\catN$ which has natural numbers $i\in\setN$ as
objects, and morphisms given by
\begin{equation*}
\catN(i,j) := 
\begin{cases}
\{\id_i\}\quad &\text{if $i = j$,}
\\
\emptyset &\text{if $i \not= j$;}
\end{cases}
\end{equation*}
it is a monoidal category with the tensor product 
\begin{equation*}
n \otimes m := n+m.
\end{equation*}
It is trivial to check that $0$ is the identity object and that
$\catN$ is freely generated by the object $1$.  Fix a vector space
\(V\). Define a functor $\EndOp[V]: \catN\opp \times \catN \to \catVect$
by
\begin{equation*}
\EndOp[V](m,n) := \Hom(V^{\otimes m},V^{\otimes n}).
\end{equation*}
This has an obvious structure of a \(\catVect\)-PROP on \(\catN\);
it is called the \emph{endomorphism PROP} of \(V\).

What should a morphism between PROPs be?  Recall that a functor $f:
\catA' \to \catA''$, i.e., a morphism of categories, is a pair of maps
$(f_{\text{Ob}}, f_{\text{Hom}})$, where $f_{\text{Hom}}$ is a natural
transform between $\catA'(-,-)$ and $\catA''(f_{\text{Ob}}(-),
f_{\text{Ob}}(-))$.  Therefore, given two PROPs ${\Pp}'$ and
${\Pp}''$, 
define a morphism $\rho: \Pp' \to \Pp''$ to be
a pair $(\rho_{\text{Ob}}, \rho_{\text{Hom}})$, where:
\begin{itemize}
\item $\rho_{\text{Ob}}: \catA'_{\text{Ob}} \to \catA''_{\text{Ob}}$
is a tensor functor,
\item $\rho_{\text{Hom}}$ is a natural transformation
$\rho_{\text{Hom}}(A,B) : {\Pp}'(A,B) \to{\Pp}''(\rho_{\text{Ob}}A,
\rho_{\text{Ob}}B)$,
\end{itemize}
that satisfies conditions that express compatibility with the tensor
structure on $\catA_{\text{Hom}}$.

If $\rho: \Pp' \to \Pp''$ is a surjective morphism, then we say that
$\Pp''$ is a PROP quotient of $\Pp'$.  For our purposes, it will
always be $\catA'_{\textrm{Ob}} = \catA''_{\textrm{Ob}}$ and
$\rho_{\textrm{Ob}} = \Id$ ; in this cases, PROP quotients are
characterized by kernels of maps $\rho_{\textrm{Hom}}(A,B): \Pp'(A,B) \to
\Pp''(A,B)$.

Now we come to the one example most relevant to this paper. Let $\Pp'
= \Pp$ be a PROP of vector spaces over $\catN$; let $\Pp''$
be the category of $\fk$-vector spaces considered as a PROP over
itself.  Since ${\catN}$ is generated by $1$ (as a monoidal category),
for any morphism
\begin{equation*}
  \rho: \Pp' \to \Pp'',  
\end{equation*}
the image of $\rho_{\textrm{Ob}}$ is generated by the vector space
$V=\rho_{\textrm{Ob}}(1)$.  Therefore, the data determining $\rho$ are
a family of morphisms
\begin{equation*}
\rho_{m,n}:{\Pp}(m,n)\to \EndOp[V](m,n).
\end{equation*}
Therefore, \(\rho\) is actually a morphism with values in the
endomorphism PROP of \(V\).
\begin{dfn}
  An \emph{action} of a linear PROP \(\Pp\) (over \(\catN\)) on a linear
  space \(V\) is a morphism \(\rho:{\Pp}\to \EndOp[V]\). The space
  \(V\) endowed with an action of \(\Pp\) is called an
  \(\Pp\)-\emph{algebra}.
\end{dfn} 
The data of an \(\Pp\)-algebra can be regarded as a family of elements
of $\Hom(V^{\otimes m},\break V^{\otimes n})$ parameterized by the space
${\Pp}(m,n)$.

\subsection{The PROP of Reshetikhin-Turaev diagrams.} 
\label{sec:RTdiagrams}

Let $L$ be the infinite strip $\setR \times [0,1]$; for $j \geq 1$ let
$s_j$, $t_j$ be the points $(j, 0) \in L$ and $(j,1) \in L$,
respectively.
\begin{figure}
\begin{equation*}
  \usegraphics{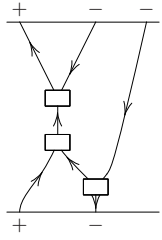}{\xy
    ,(0,-1.5);(0.6,-0.4)*+[F]{\ }**\crv{(0,-1.3)&(0.4,-0.9)%
      &(0.6,-0.4)}?(.5)*\dir{>}%
    ,(0.6,-0.4)*+[F]{\ };(0.6,0.3)*+[F]{\ }**\crv{(0.6,-0.2)}%
    ?(.75)*\dir{>}%
    ,(0.6,0.3)*+[F]{\ };(0,1.5)**\crv{(0.5,0.5)&(0.4,0.7)}%
    ?(.85)*\dir{>}%
    ,(1.2,-1.5);(1.2,-1.1)*+[F]{\ }**\crv{(1.2,-1.3)&(1.2,-1.1)}%
    ?(.10)*\dir{<}%
    ,(1.2,-1.1)*+[F]{\ };(0.6,-0.4)*+[F]{\ }**\crv{%
      (0.8,-0.7)&(0.6,-0.5)}%
    ?(.5)*\dir{>}%
    ,(0.6,0.3)*+[F]{\ };(1.2,1.5)**\crv{(0.7,0.5)&(0.8,0.7)}%
    ?(.5)*\dir{<}%
    ,(1.2,-1.1)*+[F]{\ };(2,1.5)**\crv{(1.2,-1.1)&(1.4,-0.9)&(1.5,-0.7)}%
    ?(.85)*\dir{<},%
    ,(0,-1.7)*{+},(1.2,-1.7)*{-}%
    ,(0,1.7)*{+},(1.2,1.7)*{-}%
    ,(2,1.7)*{-},%
    ,(-0.2,-1.5);(2.2,-1.5)**\dir{-}%
    ,(-0.2,1.5);(2.2,1.5)**\dir{-}%
    \endxy}
\end{equation*}
\caption{A Reshetikhin-Turaev diagram}
\end{figure}
\begin{dfn}
  A Reshetikhin-Turaev diagram\footnote{``RT-diagram'' for short.} $\Gamma$
  of type $(p,q)$ is given by a finite set $\Vertices{\Gamma}$ (the
  set of vertices) and a finite set $\Edges{\Gamma}$ (the set of
  edges) such that:
  \begin{enumerate}[RT1)]
  \item \label{RT1} each vertex \(v\) is a tiny rectangle (``coupon'' in
Reshetikhin-Turaev's original wording) contained in the strip \(L\) with
two
of its edges parallel to the boundary of \(L\) --- call them
\(\text{Top}(v)\) and \(\text{Bottom}(v)\);
\item \label{RT2} each edge is a smooth immersion \(\ell:[0,1]\to L\);
\item \label{RT3}  for each $\ell\in\Edges{\Gamma}$, the \emph{endpoints}
$\ell(0)$,
$\ell(1)$ lie in \begin{equation*}\{s_1,\ldots, s_p\}\cup\{t_1, \ldots, t_q\}
\cup\bigcup_{v\in\Vertices{\Gamma}}\bigl(\text{Top}(v)\cup\text{Bottom}(v)
\bigr);\end{equation*} the points \(\ell(0)\) and \(\ell(1)\) are called the
\emph{source} and \emph{target} of the edge \(\ell\) respectively;
\item \label{RT4} no two edges have a common endpoint;
\item \label{RT5} for every crossing of edges (including self-crossings)
an element in the set
\begin{equation*}
\left\{\, \xy*!LC\xybox{\vcross}\endxy \,,\, 
\xy*!LC\xybox{\vcrossneg}\endxy \,\right\}
\end{equation*}
must be specified, that is, we want to know which of the two
crossing arcs ``passes under''.
\end{enumerate}
\end{dfn}

A \emph{closed diagram} is a diagram of type $(0,0)$. Adding (or removing)
a  connected component of
type $(0,0)$ to a given diagram does not change its type.

\begin{xmp}
A braid on $r$ strands can be seen as a diagram of type $(r,r)$ with no
vertices (and vice-versa).
\end{xmp}

The projection $z: L \to [0,1]$ induces a differentiable function
$z\circ\ell : [0,1] \to [0,1]$ (the \emph{height function}) on every edge
$\ell\in\Edges{\Gamma}$.
\begin{dfn}
Critical points of an RT-diagram $\Gamma$ are: 
\begin{inparaenum}
\item vertices, 
\item crossings,
\item critical values of the height function on every edge.
\end{inparaenum}
\end{dfn}
If $v$ is a vertex of $\Gamma$, then let $\Leg(v)$ be the set of edges of
$\Gamma$ incident to $v$. Every set $\Leg(v)$ is divided into two
disjoint totally ordered subsets $\In(v)$ and $\Out(v)$:
\begin{equation*}
  \usegraphics{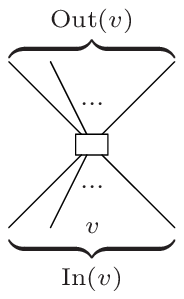}{%
    \overbrace{\underbrace{\xy%
        (0,0)*{v},
        (0,1)*+[F]{\ };%
        (-1,0)**\dir{-},(-0.5,0)**\dir{-},(0,0.5)*{\ldots},(1,0)**\dir{-},%
        (-1,2)**\dir{-},(-0.5,2)**\dir{-},(0,1.5)*{\ldots},(1,2)**\dir{-},%
        \endxy}_{\In(v)}}^{\Out(v)}
    }
\end{equation*}

A sign $\varepsilon_x \in \{ \pm1 \}$ is given to each non-critical
point $x$ according to whether the height function preserves or
reverses orientation in a neighborhood of the preimage of $x$. This
sign is locally constant (on edges except critical points); by
extension, a sign is unambiguously defined on each source and target
point; if \(p\) is an endpoint, denote its sign by \(\text{sgn}(p)\).
\begin{dfn}
The \emph{source} and the \emph{target} of an RT-diagram \(\Gamma\) of
type \((p,q)\) are
the sequences of \(\pm1\) given by 
\begin{align*}
\text{Src}(\Gamma) &:= 
(\text{sgn}(s_1),\text{sgn}(s_2),\dots,\text{sgn}(s_p)),
\text{Tgt}(\Gamma) &:= 
(\text{sgn}(t_1),\text{sgn}(t_2),\dots,\text{sgn}(t_q))
\end{align*}
For a vertex $v$ of an RT-diagram we can define $\Src(v)$ (resp.  
$\Tgt(v)$) as
the sequences of the signs of the endpoints of the edges in $\In(v)$
(resp. $\Out(v)$).
\end{dfn}

\begin{figure}[htbp]
  \begin{equation*}
    \usegraphics{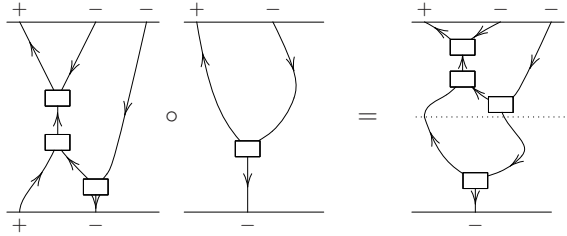}{%
      {\xy%
        ,(0,-1.5);(0.6,-0.4)*+[F]{\ }**\crv{(0,-1.3)&(0.4,-0.9)%
          &(0.6,-0.4)}?(.5)*\dir{>}%
        ,(0.6,-0.4)*+[F]{\ };(0.6,0.3)*+[F]{\ }**\crv{(0.6,-0.2)}%
        ?(.75)*\dir{>}%
        ,(0.6,0.3)*+[F]{\ };(0,1.5)**\crv{(0.5,0.5)&(0.4,0.7)}%
        ?(.85)*\dir{>}%
        ,(1.2,-1.5);(1.2,-1.1)*+[F]{\ }**\crv{(1.2,-1.3)&(1.2,-1.1)}%
        ?(.10)*\dir{<}%
        ,(1.2,-1.1)*+[F]{\ };(0.6,-0.4)*+[F]{\ }**\crv{%
          (0.8,-0.7)&(0.6,-0.5)}%
        ?(.5)*\dir{>}%
        ,(0.6,0.3)*+[F]{\ };(1.2,1.5)**\crv{(0.7,0.5)&(0.8,0.7)}%
        ?(.5)*\dir{<}%
        ,(1.2,-1.1)*+[F]{\ };(2,1.5)**\crv{(1.2,-1.1)&(1.4,-0.9)&(1.5,-0.7)}%
        ?(.85)*\dir{<},%
        ,(0,-1.7)*{+},(1.2,-1.7)*{-}%
        ,(0,1.7)*{+},(1.2,1.7)*{-}%
        ,(2,1.7)*{-},%
        ,(-0.2,-1.5);(2.2,-1.5)**\dir{-}%
        ,(-0.2,1.5);(2.2,1.5)**\dir{-}%
        \endxy}\space\circ\space
      {\xy%
        ,(0.8,-1.5);(0.8,-0.5)*+[F]{\ }**\dir{-}
        ?(.4)*\dir{<}%
        ,(0.8,-0.5)*+[F]{\ };(0,1.5)**\crv{(0.4,-0.1)&(0.2,0.3)}%
        ?(.85)*\dir{>}%
        ,(0.8,-0.5)*+[F]{\ };(1.2,1.5)**\crv{(0.9,-0.4)&%
          (1.2,-0.3)&(1.6,0.3)&(1.6,0.7)}%
        ?(.85)*\dir{<}%
        ,(0.8,-1.7)*{-},(0,1.7)*{+}%
        ,(1.2,1.7)*{-}%
        ,(-0.2,-1.5);(2,-1.5)**\dir{-}%
        ,(-0.2,1.5);(2,1.5)**\dir{-}%
        \endxy}
      \quad = \quad
      {\xy%
        ,(0.8,-1.5);(0.8,-1)*+[F]{\ }**\dir{-}%
        ?(.25)*\dir{<}%
        ,(0.8,-1)*+[F]{\ };(0,0)**\crv{(0.4,-0.8)&(0.2,-0.6)}%
        ?(.85)*\dir{>}%
        ,(0.8,-1)*+[F]{\ };(1.2,0.2)*+[F]{\ }**\crv{(0.9,-0.95)&%
          (1.2,-0.9)&(1.6,-0.6)&(1.6,-0.4)&(1.2,0.0)}%
        ?(.35)*\dir{<}%
        ,(0,0);(0.6,0.6)*+[F]{\ }**\crv{(0,0.2)&(0.4,0.3)%
          &(0.6,0.6)}%
        ,(0.6,0.6)*+[F]{\ };(0.6,1.1)*+[F]{\ }**\dir{-}%
        ?(.75)*\dir{>}%
        ,(0.6,1.1)*+[F]{\ };(0,1.5)**\crv{(0.5,1.1)&(0.4,1.3)}%
        ?(.85)*\dir{>}%
        ,(1.2,0.2)*+[F]{\ };(0.6,0.6)*+[F]{\ }**\crv{%
          (0.8,0.3)&(0.6,0.7)}%
        ?(.5)*\dir{>}%
        ,(0.6,1.1)*+[F]{\ };(1.2,1.5)**\crv{(0.7,1.2)&(0.8,1.3)}%
        ?(.5)*\dir{<}%
        ,(1.2,0.2)*+[F]{\ };(2,1.5)**\crv{(1.2,0.2)&(1.4,0.3)&(1.5,0.4)}%
        ?(.85)*\dir{<},%
        ,(-0.2,0);(2.2,0)**\dir{.}%
        ,(0,1.7)*{+},(1.2,1.7)*{-}%
        ,(2,1.7)*{-},%
        ,(-0.2,1.5);(2.2,1.5)**\dir{-}%
        ,(0.8,-1.7)*{-}%
        ,(-0.2,-1.5);(2,-1.5)**\dir{-}%
        \endxy}
      }
  \end{equation*}
\caption{Composition product of graphs.}
\label{fig:graph-composition}
\end{figure}

\begin{dfn}\label{dfn:graph-composition}
Two RT-diagrams $\Gamma$ and $\Phi$ are \emph{composable} iff
$\Src(\Gamma) = \Tgt(\Phi)$. If $\Gamma$ and $\Phi$ are composable, we
can
form a new diagram $\Gamma\circ\Phi$ by ``stacking $\Gamma$ on
top of $\Phi$''
(see \prettyref{fig:graph-composition}).
\end{dfn}

Note that this composition product restricts to the usual braid
composition on diagrams corresponding to elements of the braid group.

\begin{figure}[htbp]
  \begin{equation*}
    \usegraphics{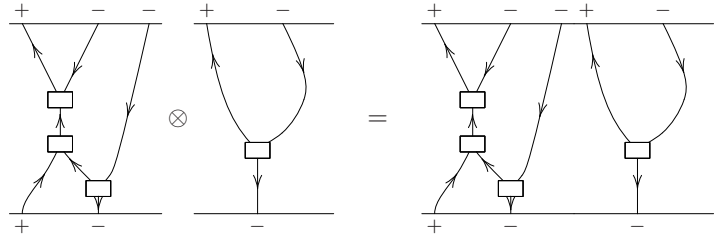}{%
      {\xy%
        ,(0,-1.5);(0.6,-0.4)*+[F]{\ }**\crv{(0,-1.3)&(0.4,-0.9)%
          &(0.6,-0.4)}?(.5)*\dir{>}%
        ,(0.6,-0.4)*+[F]{\ };(0.6,0.3)*+[F]{\ }**\crv{(0.6,-0.2)}%
        ?(.75)*\dir{>}%
        ,(0.6,0.3)*+[F]{\ };(0,1.5)**\crv{(0.5,0.5)&(0.4,0.7)}%
        ?(.85)*\dir{>}%
        ,(1.2,-1.5);(1.2,-1.1)*+[F]{\ }**\crv{(1.2,-1.3)&(1.2,-1.1)}%
        ?(.10)*\dir{<}%
        ,(1.2,-1.1)*+[F]{\ };(0.6,-0.4)*+[F]{\ }**\crv{%
          (0.8,-0.7)&(0.6,-0.5)}%
        ?(.5)*\dir{>}%
        ,(0.6,0.3)*+[F]{\ };(1.2,1.5)**\crv{(0.7,0.5)&(0.8,0.7)}%
        ?(.5)*\dir{<}%
        ,(1.2,-1.1)*+[F]{\ };(2,1.5)**\crv{(1.2,-1.1)&(1.4,-0.9)&(1.5,-0.7)}%
        ?(.85)*\dir{<},%
        ,(0,-1.7)*{+},(1.2,-1.7)*{-}%
        ,(0,1.7)*{+},(1.2,1.7)*{-}%
        ,(2,1.7)*{-},%
        ,(-0.2,-1.5);(2.2,-1.5)**\dir{-}%
        ,(-0.2,1.5);(2.2,1.5)**\dir{-}%
        \endxy}\space\otimes\space
      {\xy%
        ,(0.8,-1.5);(0.8,-0.5)*+[F]{\ }**\dir{-}
        ?(.4)*\dir{<}%
        ,(0.8,-0.5)*+[F]{\ };(0,1.5)**\crv{(0.4,-0.1)&(0.2,0.3)}%
        ?(.85)*\dir{>}%
        ,(0.8,-0.5)*+[F]{\ };(1.2,1.5)**\crv{(0.9,-0.4)&%
          (1.2,-0.3)&(1.6,0.3)&(1.6,0.7)}%
        ?(.85)*\dir{<}%
        ,(0.8,-1.7)*{-},(0,1.7)*{+}%
        ,(1.2,1.7)*{-}%
        ,(-0.2,-1.5);(2,-1.5)**\dir{-}%
        ,(-0.2,1.5);(2,1.5)**\dir{-}%
        \endxy}
      \quad = \quad
      {\xy%
        ,(0,-1.5);(0.6,-0.4)*+[F]{\ }**\crv{(0,-1.3)&(0.4,-0.9)%
          &(0.6,-0.4)}?(.5)*\dir{>}%
        ,(0.6,-0.4)*+[F]{\ };(0.6,0.3)*+[F]{\ }**\crv{(0.6,-0.2)}%
        ?(.75)*\dir{>}%
        ,(0.6,0.3)*+[F]{\ };(0,1.5)**\crv{(0.5,0.5)&(0.4,0.7)}%
        ?(.85)*\dir{>}%
        ,(1.2,-1.5);(1.2,-1.1)*+[F]{\ }**\crv{(1.2,-1.3)&(1.2,-1.1)}%
        ?(.10)*\dir{<}%
        ,(1.2,-1.1)*+[F]{\ };(0.6,-0.4)*+[F]{\ }**\crv{%
          (0.8,-0.7)&(0.6,-0.5)}%
        ?(.5)*\dir{>}%
        ,(0.6,0.3)*+[F]{\ };(1.2,1.5)**\crv{(0.7,0.5)&(0.8,0.7)}%
        ?(.5)*\dir{<}%
        ,(1.2,-1.1)*+[F]{\ };(2,1.5)**\crv{(1.2,-1.1)&(1.4,-0.9)&(1.5,-0.7)}%
        ?(.85)*\dir{<},%
        ,(0,-1.7)*{+},(1.2,-1.7)*{-}%
        ,(0,1.7)*{+},(1.2,1.7)*{-}%
        ,(2,1.7)*{-},%
        ,(-0.2,-1.5);(2.2,-1.5)**\dir{-}%
        ,(-0.2,1.5);(2.2,1.5)**\dir{-}%
        \endxy}
      {\xy%
        ,(0.8,-1.5);(0.8,-0.5)*+[F]{\ }**\dir{-}
        ?(.4)*\dir{<}%
        ,(0.8,-0.5)*+[F]{\ };(0,1.5)**\crv{(0.4,-0.1)&(0.2,0.3)}%
        ?(.85)*\dir{>}%
        ,(0.8,-0.5)*+[F]{\ };(1.2,1.5)**\crv{(0.9,-0.4)&%
          (1.2,-0.3)&(1.6,0.3)&(1.6,0.7)}%
        ?(.85)*\dir{<}%
        ,(0.8,-1.7)*{-},(0,1.7)*{+}%
        ,(1.2,1.7)*{-}%
        ,(-0.2,-1.5);(2,-1.5)**\dir{-}%
        ,(-0.2,1.5);(2,1.5)**\dir{-}%
        \endxy}
      }
  \end{equation*}
\caption{Tensor product of graphs.}
\label{fig:graph-otimes}
\end{figure}

By the above definition, we can take the linear spans $\Dia(S,T)$ of
diagrams with given source $S$ and target $T$, as the $\Hom$-spaces of a
suitable PROP.
\begin{dfn}
\label{dfn:graph-category}
$\Dia$ is the PROP which has finite sequences of $\pm1$ as objects; the
$\Hom$-space $\Dia(S,T)$ is the linear span of the set of diagrams
with source $S$ and target $T$.

Composition of morphisms is defined by bilinear extension of the
composition product $\circ$ (see \prettyref{dfn:graph-composition}).

The tensor product is given on objects by concatenation:
\begin{equation*}
(\varepsilon_1, \ldots, \varepsilon_r) \otimes (\varepsilon'_1,
\ldots, \varepsilon'_s) = (\varepsilon_1, \ldots, \varepsilon_r,
\varepsilon'_1,
\ldots, \varepsilon'_s),
\end{equation*}
and on morphisms by juxtaposition of diagrams (see
\prettyref{fig:graph-otimes}). 

The braiding is given by diagrams corresponding to elements of the braid
groups. 
\end{dfn}

Any RT-diagram can be considered as the planar projection of a purely
1-dimensional CW-complex with oriented edges, and a partition of
half-edges occurring at each vertex and of its endpoints into two
totally ordered disjoint subsets.  The class of CW-complexes with such
an additional structure is the class of \emph{RT-graphs}. There is a
natural map ``forgetting the planar immersion'' $\Dia\to\RT$. This
forgetful functor is actually a PROP quotient, as the following lemmas
state.
\begin{lemma}[\cite{reshetikhin-turaev;ribbon-graphs}]
\label{lemma:generators}
The PROP \(\Dia\) of Reshetikhin-Turaev diagrams is the free PROP
generated by the following elementary pieces (note that an orientation
must be added to the strands!)
\begin{center}
  \usegraphics{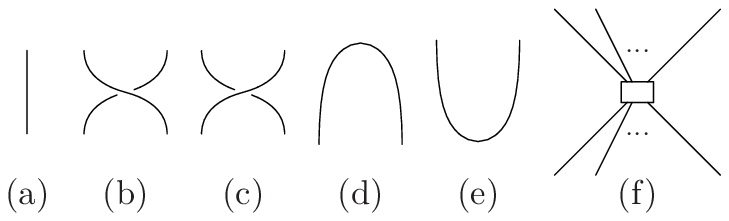}{%
    \begin{tabular}{cccccc}
      \(\xy*!LC\xybox{(0,0);(0,1)**\dir{-}}\endxy\)
      &
      \(\xy*!LC\xybox{%
        \vcross~{(0,1)}{(1,1)}{(0,0)}{(1,0)}}\endxy\)
      &
      \(\xy*!LC\xybox{%
        \vcross~{(0,0)}{(1,0)}{(0,1)}{(1,1)}}\endxy\)
      &
      \(\xy*!LC\xybox{%
        \vloop~{(0,1)}{(1,1)}{(0,0)}{(1,0)}}\endxy\)
      &
      \(\xy*!LC\xybox{%
        \vloop~{(0,0)}{(1,0)}{(0,1)}{(1,1)}}\endxy\)
      &
      \(\xy*!LC\xybox{
        (0,1)*+[F]{\ };%
        (-1,0)**\dir{-},(-0.5,0)**\dir{-},%
        (0,0.5)*+{\ldots},(1,0)**\dir{-},%
        (-1,2)**\dir{-},(-0.5,2)**\dir{-},%
        (0,1.5)*+{\ldots},(1,2)**\dir{-},%
        }\endxy\)
      \\
      (a) & (b) & (c) & (d) & (e) & (f)
    \end{tabular}
    }
\end{center}
\end{lemma}
A piece of type (a) is called a ``strand''; those of type (b) and (c)
are named ``crossings''; (d) and (e) are the ``coupling'' and the
``Casimir''; (f) is, plainly, a ``vertex''.
\begin{rem} \prettyref{lemma:generators} just states that an
RT-diagram is a composition of ``rows'' made of pieces
of type (a)--(f). Generically, such rows will be made of one piece of
type (b)--(f) padded with a number of strands (a)
on the two sides.
\end{rem}

\begin{figure}[htbp]
  \begin{center}
    \usegraphics{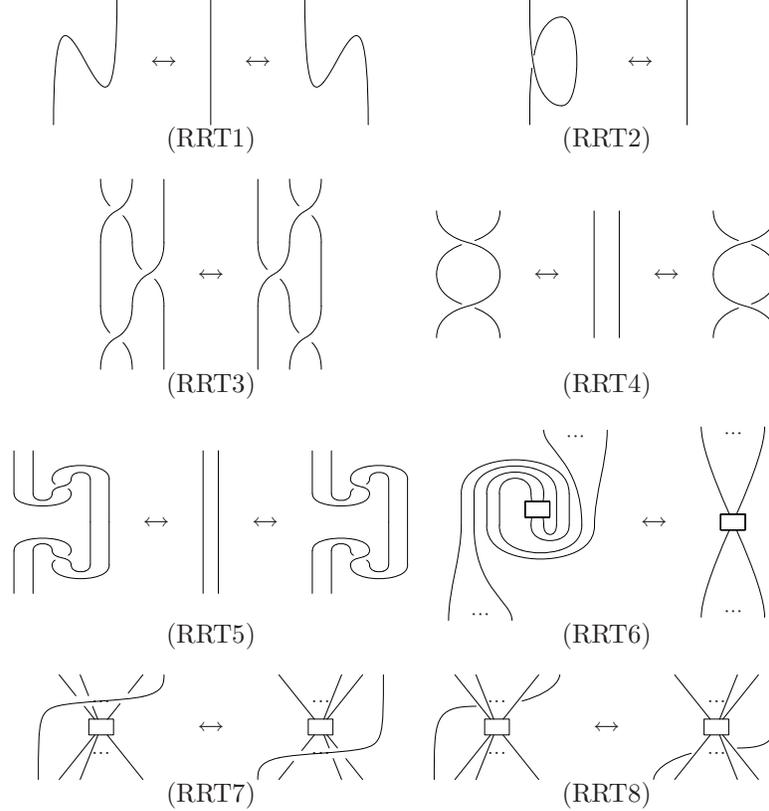}{%
      \begin{tabular}{cc}
        \(
        {\xy(-.5,-1);(.5,1)**\crv{(-.5,.5)&(-.2,0.5)&(.2,-0.5)&(.5,-.5)}\endxy}
        \quad\leftrightarrow\quad
        {\xy(0,-1);(0,1)**\dir{-}\endxy}
        \quad\leftrightarrow\quad
        {\xy(.5,-1);(-.5,1)**\crv{(.5,.5)&(.2,0.5)&(-.2,-0.5)&(-.5,-.5)}\endxy}
        \)
        &
        \(
        {\xy
          \hunder~{(0,.7)}{(.5,.7)}{(0,0.2)}{(.5,.2)}
          \hcap~{(0.5,.7)}{(1,.7)}{(0.5,0.2)}{(1,0.2)}
          \hover~{(0,-.2)}{(.5,-.2)}{(0,-.7)}{(.5,-.7)}
          \hcap~{(0.5,-.2)}{(1,-.2)}{(0.5,-.7)}{(1,-.7)}
          ,(0,.7);(0,1)**\dir{-},(0,.2);(0,-.2)**\dir{-},(0,-.7);(0,-1)**\dir{-}
          \endxy}
        \quad\leftrightarrow\quad
        {\xy
          ,(0,-1);(0,1)**\dir{-}
          \endxy}
        \)
        \\
        (RRT1)
        &
        (RRT2)
        \\[2ex]
        \(
        {\xy
          \vcrossneg~{(0,1.5)}{(.5,1.5)}{(0,0.5)}{(.5,0.5)}
          \vcrossneg~{(.5,0.5)}{(1,0.5)}{(.5,-.5)}{(1,-.5)}
          \vcrossneg~{(0,-0.5)}{(.5,-0.5)}{(0,-1.5)}{(.5,-1.5)}
          ,(1,1.5);(1,.5)**\dir{-},(0,.5);(0,-.5)**\dir{-},
          (1,-.5);(1,-1.5)**\dir{-}
          \endxy}
        \quad\leftrightarrow\quad
        {\xy
          \vcrossneg~{(0.5,1.5)}{(1,1.5)}{(0.5,0.5)}{(1,0.5)}
          \vcrossneg~{(0,0.5)}{(.5,0.5)}{(0,-.5)}{(.5,-.5)}
          \vcrossneg~{(0.5,-0.5)}{(1,-0.5)}{(0.5,-1.5)}{(1,-1.5)}
          ,(0,1.5);(0,.5)**\dir{-},(1,.5);(1,-.5)**\dir{-},
          (0,-.5);(0,-1.5)**\dir{-}
          \endxy}
        \)
        &
        \(
        {\xy*!LC\xybox{\vcross\vcrossneg}\endxy}
        \quad\leftrightarrow\quad
        {\xy(-.2,-1);(-.2,1)**\dir{-},(.2,-1);(.2,1)**\dir{-}\endxy}
        \quad\leftrightarrow\quad
        {\xy*!LC\xybox{\vcrossneg\vcross}\endxy}
        \)
        \\
        (RRT3)
        &
        (RRT4)
        \\[2ex]
        \({\xy(0.6,0):(0,0.75)::
          \vcrossneg~{(0,1.5)}{(.5,1.5)}{(0,1)}{(.5,1)}
          \vcap~{(0.5,2)}{(1,2)}{(0.5,1.5)}{(1,1.5)}
          \vcap~{(-0.5,0.5)}{(0,0.5)}{(-0.5,1)}{(0,1)}
          \vcap~{(0,2.5)}{(1.5,2.5)}{(0,1.5)}{(1.5,1.5)}
          \vcap~{(-1,0.1)}{(0.5,0.1)}{(-1,1)}{(0.5,1)}
          \vcrossneg~{(0,-1.5)}{(.5,-1.5)}{(0,-1)}{(.5,-1)}
          \vcap~{(0.5,-2)}{(1,-2)}{(0.5,-1.5)}{(1,-1.5)}  
          \vcap~{(-0.5,-.5)}{(0,-.5)}{(-0.5,-1)}{(0,-1)}
          \vcap~{(0,-2.5)}{(1.5,-2.5)}{(0,-1.5)}{(1.5,-1.5)}     
          \vcap~{(-1,-0.1)}{(0.5,-0.1)}{(-1,-1)}{(0.5,-1)},
          (-1,1);(-1,2.5)**\dir{-},(1,1.5);(1,0)**\dir{-},
          (1.5,1.5);(1.5,0)**\dir{-},
          (-1,-1);(-1,-2.5)**\dir{-},(1,-1.5);(1,0)**\dir{-},
          (1.5,-1.5);(1.5,0)**\dir{-},
          (-.5,1);(-.5,2.5)**\dir{-},(-.5,-1);(-.5,-2.5)**\dir{-}
          \endxy}
        \quad\leftrightarrow\quad
        {\xy(0.6,0):(0,0.75)::
          (-0.2,2.5);(-0.2,-2.5)**\dir{-},
          (0.2,2.5);(0.2,-2.5)**\dir{-}\endxy}
        \quad\leftrightarrow\quad
        {\xy(0.6,0):(0,0.75)::
          \vcross~{(0,1.5)}{(.5,1.5)}{(0,1)}{(.5,1)}
          \vcap~{(0.5,2)}{(1,2)}{(0.5,1.5)}{(1,1.5)}
          \vcap~{(-0.5,0.5)}{(0,0.5)}{(-0.5,1)}{(0,1)}
          \vcap~{(0,2.5)}{(1.5,2.5)}{(0,1.5)}{(1.5,1.5)}
          \vcap~{(-1,0.1)}{(0.5,0.1)}{(-1,1)}{(0.5,1)}
          \vcross~{(0,-1.5)}{(.5,-1.5)}{(0,-1)}{(.5,-1)}
          \vcap~{(0.5,-2)}{(1,-2)}{(0.5,-1.5)}{(1,-1.5)}  
          \vcap~{(-0.5,-.5)}{(0,-.5)}{(-0.5,-1)}{(0,-1)}
          \vcap~{(0,-2.5)}{(1.5,-2.5)}{(0,-1.5)}{(1.5,-1.5)}     
          \vcap~{(-1,-0.1)}{(0.5,-0.1)}{(-1,-1)}{(0.5,-1)},
          (-1,1);(-1,2.5)**\dir{-},(1,1.5);(1,0)**\dir{-},
          (1.5,1.5);(1.5,0)**\dir{-},
          (-1,-1);(-1,-2.5)**\dir{-},(1,-1.5);(1,0)**\dir{-},
          (1.5,-1.5);(1.5,0)**\dir{-},
          (-.5,1);(-.5,2.5)**\dir{-},(-.5,-1);(-.5,-2.5)**\dir{-}
          \endxy}
        \)
        &
        \(
        {\xy*!LC\xybox{
            \vcap~{(-0.5,1)}{(0,1)}{(-0.5,0.5)}{(0,0.5)}
            \vcap~{(-0.7,1.2)}{(0.2,1.2)}{(-0.7,0.5)}{(0.2,0.5)}
            \vcap~{(-0.9,1.4)}{(0.4,1.4)}{(-0.9,0.5)}{(0.4,0.5)}
            \vcap~{(-1.1,1.6)}{(0.6,1.6)}{(-1.1,0.5)}{(0.6,0.5)}
            \vcap~{(0,-.5)}{(0.6,-.5)}{(0,0)}{(0.6,0)}  
            \vcap~{(0.2,-.3)}{(0.4,-.3)}{(0.2,0)}{(0.4,0)}  
            \vcap~{(-0.5,-.9)}{(0.8,-.9)}{(-0.5,0)}{(0.8,0)}  
            \vcap~{(-0.7,-1.1)}{(1,-1.1)}{(-0.7,0)}{(1,0)}
            ,(0.1,0.25)*+[F]{\ };(0,0.5)**\crv{(0,0.3)}
            ,(0.1,0.25)*+[F]{\ };(0.2,0.5)**\crv{(0.2,0.3)}
            ,(0.1,0.25)*+[F]{\ };(0,0)**\crv{(0,0.2)}
            ,(0.1,0.25)*+[F]{\ };(0.2,0)**\crv{(0.2,0.2)}
            ,(0.8,0);(0.2,1.5)**\crv{(0.8,1)&(0.2,1.2)}
            ,(1,0);(1.2,1.5)**\crv{(1,0.5)&(1.2,1)}
            ,(0.7,1,4)*{\cdots}
            ,(-0.9,0);(-0.3,-1.5)**\crv{(-0.9,-1)&(-0.3,-1.2)}
            ,(-1.1,0);(-1.3,-1.5)**\crv{(-1.1,-0.5)&(-1.3,-1)}
            ,(-0.8,-1,4)*{\cdots}
            ,(0.4,0.5);(0.4,0)**\dir{-}
            ,(0.6,0.5);(0.6,0)**\dir{-}
            ,(-0.5,0.5);(-0.5,0)**\dir{-}
            ,(-0.7,0.5);(-0.7,0)**\dir{-}
            ,(-0.9,0.5);(-0.9,0)**\dir{-}
            ,(-1.1,0.5);(-1.1,0)**\dir{-}
            }
          \endxy}
        \quad\leftrightarrow\quad
        {\xy*!LC\xybox{
            ,(0,0)*+[F]{\ };(-0.5,-1.5)**\crv{(-0.5,-1.1)}
            ,(0,0)*+[F]{\ };(0.5,-1.5)**\crv{(0.5,-1.1)}
            ,(0,0)*+[F]{\ };(-0.5,1.5)**\crv{(-0.5,1.1)}
            ,(0,0)*+[F]{\ };(0.5,1.5)**\crv{(0.5,1.1)}
            ,(0,1.4)*{\cdots},(0,-1.4)*{\cdots}
            }
          \endxy}
        \)
        \\
        (RRT5)
        &
        (RRT6)
        \\[2ex]
        \(
        {\xy
          *!LC\xybox{(0.66,0):(0,1.25)::
            (0,1)*+[F]{\ };%
            (-1,0)**\dir{-},(-0.5,0)**\dir{-},%
            (0,0.5)*+{\ldots},(1,0)**\dir{-},%
            (-.4,1.4)**\dir{-},(-0.2,1.4)**\dir{-},%
            (0,1.5)*+{\ldots},(0.45,1.45)**\dir{-},%
            ,(-1.5,0);(1.5,2)**\crv{(-1.5,1.15)&(-1.5,1.5)&(1.5,1.5)&(1.5,1.85)}
            ,(-.6,1.6);(-1,2)**\dir{-},(-0.3,1.6);(-0.5,2)**\dir{-}
            ,(0.65,1.65);(1,2)**\dir{-}
            }
          \endxy}
        \quad\leftrightarrow\quad
        {\xy
          *!LC\xybox{(0.66,0):(0,1.25)::
            (0,-1)*+[F]{\ };%
            (1,0)**\dir{-},(0.5,0)**\dir{-},%
            (0,-0.5)*+{\ldots},(-1,0)**\dir{-},%
            (.4,-1.4)**\dir{-},(0.2,-1.4)**\dir{-},%
            (0,-1.5)*+{\ldots},(-0.45,-1.45)**\dir{-},%
            ,(1.5,0);(-1.5,-2)**\crv{(1.5,-1.15)&(1.5,-1.5)%
              &(-1.5,-1.5)&(-1.5,-1.85)}
            ,(.6,-1.6);(1,-2)**\dir{-},(0.3,-1.6);(0.5,-2)**\dir{-}
            ,(-0.65,-1.65);(-1,-2)**\dir{-}
            }
          \endxy}
        \)
        &
        \( 
        {\xy
          *!LC\xybox{(0.66,0):(0,1.25)::
            (0,1)*+[F]{\ };%
            (-1,0)**\dir{-},(-0.5,0)**\dir{-},%
            (0,0.5)*+{\ldots},(1,0)**\dir{-},%
            (-1,2)**\dir{-},(-0.5,2)**\dir{-},%
            (0,1.5)*+{\ldots},(1,2)**\dir{-},%
            ,(-1.5,0);(-.5,1.4)**\crv{(-1.5,1.15)&(-1.5,1.4)}
            ,(0.6,1.5);(1.5,2)**\crv{(1.5,1.7)&(1.5,2)}
            }
          \endxy}
        \quad\leftrightarrow\quad
        {\xy
          *!LC\xybox{(0.66,0):(0,1.25)::
            (0,-1)*+[F]{\ };%
            (1,0)**\dir{-},(0.5,0)**\dir{-},%
            (0,-0.5)*+{\ldots},(-1,0)**\dir{-},%
            (1,-2)**\dir{-},(0.5,-2)**\dir{-},%
            (0,-1.5)*+{\ldots},(-1,-2)**\dir{-},%
            ,(1.5,0);(.5,-1.4)**\crv{(1.5,-1.15)&(1.5,-1.4)}
            ,(-0.6,-1.5);(-1.5,-2)**\crv{(-1.5,-1.7)&(-1.5,-2)}
            }
          \endxy}
        \)
        \\
        (RRT7)
        &
        (RRT8)
      \end{tabular}
      }
    \caption{Reidmeister-Reshetikhin-Turaev moves for
      RT-diagrams. Move (RRT2) is different if the ground category has
      non-trivial balancing.}
    \label{fig:rrt}
  \end{center}
\end{figure}

\begin{lemma}[\cite{reshetikhin-turaev;ribbon-graphs}]\label{lemma:moves}
  The PROP \(\RT\) of RT-graphs is the quotient of \(\Dia\) with
  respect to the following relations:
\begin{enumerate}[M1)]
\item \label{M1}a finite set of graphical moves, called the
  Reidmeister-Reshetikhin-Turaev moves (see \prettyref{fig:rrt} or
  \cite{reshetikhin-turaev;ribbon-graphs} for a listing);
\item \label{M2} the equivalence of undercrossings and overcrossings.
\end{enumerate}
\end{lemma}
\begin{rem}
  Notice that the list of Reidmeister-Reshetikhin-Turaev moves is
  different if the ground category has non-trivial balancing
  (cf. \cite{joyal-street;1991} and \cite{shum;tortile-categories}),
  indeed one would need graphs with ``framed'' edges, and the move (RRT2)
  introduces a twist.
\end{rem}
\begin{rem}
\prettyref{lemma:moves} states that if \(\Gamma_1\) and \(\Gamma_2\) are
two planar projections that realize an RT-graph \(\Gamma\) as an
RT-diagram,
then one can change \(\Gamma_1\) into \(\Gamma_2\) by a finite sequence of
moves M\ref{M1}-M\ref{M2}.
\end{rem}

\subsection{Graphical
  calculus on rigid braided tensor categories} \label{sec:rt-graphs}
Now let $\catA$ be a tensor category. Define $\freemsc{\catA}$ to be
the category whose objects are finite sequences $(A_1, \ldots, A_r;
\varepsilon_1, \ldots, \varepsilon_r)$ of objects in $\catA$ and signs
$\pm1$, whereas a morphism $(A_*, \varepsilon_*) \to (B_*, \delta_*)$
is an element $f \in \catA(A_1^{\varepsilon_1} \otimes \ldots \otimes
A_r^{\varepsilon_r}, B_1^{\delta_1} \otimes \ldots B_s^{\delta_s})$
--- recall that $A^1 = A$ and $A^{-1} = \rdual{A}$.  There is an
obvious functor $\freemsc{\catA} \to \catA$ defined on objects by
$(A_1, \ldots, A_r; \varepsilon_1, \ldots, \varepsilon_r) \mapsto
A_1^{\varepsilon_1} \otimes \cdots \otimes A_r^{\varepsilon_r}$.
\begin{dfn}
An $\catA$-colored RT-diagram $\Gamma$ is an RT-diagram together with
\begin{enumerate}
\item an assignment of an object $A_\ell \in \catA$ for each
$\ell\in\Edges{\Gamma}$;
\item an assignment of a morphism $f\in\catA(\Src(v),\Tgt(v))$ for each
vertex $v$ of $\Gamma$, where
\(\Src(v)\) and \(\Tgt(v)\) are the sequences
\((A_1,\dots,A_r;\varepsilon_1, \ldots,
\varepsilon_r)\) of objects and signs decorating edges in
\(\In(v)\) and \(\Out(v)\).
\end{enumerate}
The source and the target of an $\catA$-colored RT-diagram are defined
analogously
and are denoted $\Src(\Gamma)$ and $\Tgt(\Gamma)$
respectively.
\end{dfn}
It is trivial to generalize \prettyref{dfn:graph-category} to
\(\catA\)-colored RT-diagrams; call
$\Dia[\catA]$ the PROP of $\catA$-colored RT
diagrams. $\Src{}_{\catA}(\Gamma)$
and $\Tgt{}_{\catA}(\Gamma)$
are objects of $\freemsc{\catA}$ for each $\catA$-colored RT-diagram.
Obviously, $\circ$ and $\otimes$ are bilinear with respect to the vector
space structure on $\Dia[\catA]$. Any \(\catA\)-colored RT-diagram
can be seen as the planar projection of an \(\catA\)-colored RT-graph
and the following analogue of lemmas \ref{lemma:generators} and
\ref{lemma:moves}
hold.
\begin{lemma}
\label{lemma:A-generators}
The PROP $\Dia[\catA]$ of $\catA$-colored RT-diagrams is the free
PROP
generated by the following elementary pieces (note that an orientation
must be added to the strands!)
\begin{center}
  \usegraphics{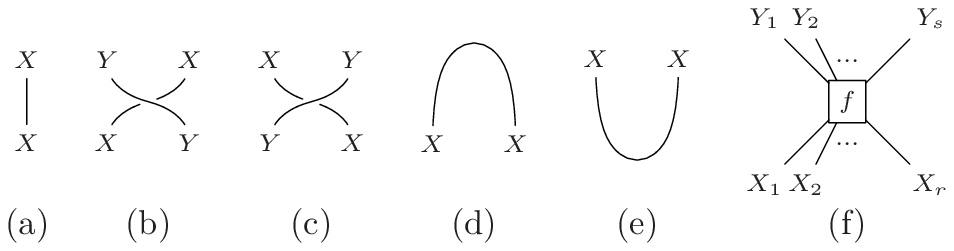}{%
    \begin{tabular}{cccccc}
      \(\xy*!LC\xybox{(0,0)*+{X};(0,1)*+{X}**\dir{-}}\endxy\)
      &
      \(\xy*!LC\xybox{%
        \vcross~{(0,1)*+{Y}}{(1,1)*+{X}}{(0,0)*+{X}}{(1,0)*+{Y}}}\endxy\)
      &
      \(\xy*!LC\xybox{%
        \vcross~{(0,0)*+{Y}}{(1,0)*+{X}}{(0,1)*+{X}}{(1,1)*+{Y}}}\endxy\)
      &
      \(\xy*!LC\xybox{%
        \vloop~{(0,1)}{(1,1)}{(0,0)*+{X}}{(1,0)*+{X}}}\endxy\)
      &
      \(\xy*!LC\xybox{%
        \vloop~{(0,0)}{(1,0)}{(0,1)*+{X}}{(1,1)*+{X}}}\endxy\)
      &
      \(\xy*!LC\xybox{
        (0,1)*+[F]{f};%
        (-1,0)*+{X_1}**\dir{-},(-0.5,0)*+{X_2}**\dir{-},%
        (0,0.5)*+{\ldots},(1,0)*+{X_r}**\dir{-},%
        (-1,2)*+{Y_1}**\dir{-},(-0.5,2)*+{Y_2}**\dir{-},%
        (0,1.5)*+{\ldots},(1,2)*+{Y_s}**\dir{-},%
        }\endxy\)
      \\
      (a) & (b) & (c) & (d) & (e) & (f)
    \end{tabular}
    }
\end{center}
\end{lemma}
\begin{lemma}\label{lemma:Amoves}
The PROP \(\RT[\catA]\) of \(\catA\)-colored RT-graphs is the quotient
of \(\Dia[\catA]\) with respect to the
relations generated by moves M\ref{M1}--M\ref{M2}.
\end{lemma}

Since \(\Dia[\catA]\) is free as a PROP, to give a tensor functor
\(\Dia[\catA]\to{\category{B}}\) it suffices to define it on the
generators. In particular, taking \(\category{B}=\catA\) we find the
following
\begin{prop}[Reshetikhin-Turaev's graphical calculus
\cite{reshetikhin-turaev;ribbon-graphs}]
For any rigid braided tensor category \(\catA\), there is a tensor functor
\(Z_{\catA}:
\Dia[\catA] \to \catA\), mapping an object \((A_1, \ldots, A_r;
\varepsilon_1,
\ldots, \varepsilon_k) \in \Dia[\catA]\) to \(A_1^{\varepsilon_1}
\otimes
\dots \otimes A_k^{\varepsilon_k}
\in \catA\), and defined on generators of morphisms in \(\Dia[\catA]\)
as
\begin{center}
  \everyxy={/r24pt/:}
  \usegraphics{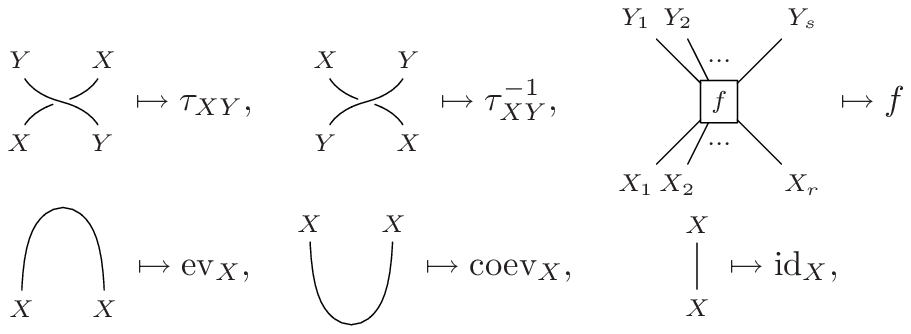}{%
    \begin{tabular}{ccc}
      \(\xy*!LC\xybox{%
        \vcross~{(0,1)*+{Y}}{(1,1)*+{X}}{(0,0)*+{X}}{(1,0)*+{Y}}}\endxy
      \mapsto \tau_{XY},\)
      \label{graph-cross+}
      &
      \(\xy*!LC\xybox{%
        \vcross~{(0,0)*+{Y}}{(1,0)*+{X}}{(0,1)*+{X}}{(1,1)*+{Y}}}\endxy
      \mapsto \tau_{XY}^{-1},\)
      \label{graph-cross-}
      &
      \(\xy*!LC\xybox{
        (0,1)*+[F]{f};%
        (-1,0)*+{X_1}**\dir{-},(-0.5,0)*+{X_2}**\dir{-},%
        (0,0.5)*+{\ldots},(1,0)*+{X_r}**\dir{-},%
        (-1,2)*+{Y_1}**\dir{-},(-0.5,2)*+{Y_2}**\dir{-},%
        (0,1.5)*+{\ldots},(1,2)*+{Y_s}**\dir{-},%
        }\endxy \label{graph-morphism} \mapsto f\)
      \\
      \(\xy*!LC\xybox{%
        \vloop~{(0,1)}{(1,1)}{(0,0)*+{X}}{(1,0)*+{X}}}\endxy \mapsto
      \ev_{X},\)
      \label{graph-casimir}
      &
      \(\xy*!LC\xybox{%
        \vloop~{(0,0)}{(1,0)}{(0,1)*+{X}}{(1,1)*+{X}}}\endxy \mapsto
      \coev_{X},\)
      \label{graph-coupling}
      &
      \(\xy*!LC\xybox{(0,0)*+{X};(0,1)*+{X}**\dir{-}}\endxy \mapsto
      \id_X,\)
      \label{graph-id}
    \end{tabular}
    }
  \end{center}
where \(\tau_{XY}\), \(\ev_X\), \(\coev_X\) are the structure maps of
\(\catA\), and \(f\) is a morphism in \(\catA\); take the dual of an
object if the sign $\varepsilon$ on the corresponding edge is $-1$.
\end{prop}
\begin{rem}
By definition of a tensor functor, the following relations hold:
\begin{equation*}
Z_{\catA}(\Gamma\circ\Phi) = Z_{\catA}(\Gamma) \circ
Z_{\catA}(\Phi), \qquad
Z_{\catA}(\Gamma\otimes\Phi) = Z_{\catA}(\Gamma) \otimes
Z_{\catA}(\Phi).
\end{equation*}
Moreover, $Z_{\catA}$ is linear:
\begin{equation*}
Z_{\catA}(a\Gamma + b\Phi) = aZ_{\catA}(\Gamma) + bZ_{\catA}(\Phi).
\end{equation*}
\end{rem}
Now we come to the main result of this section.
\begin{thm}[Reshetikhin-Turaev,
\cite{reshetikhin-turaev;ribbon-graphs}]
\label{thm:rt1}
Let \(\catA\) be a rigid symmetric tensor category. Then the
Reshetikhin-Turaev's
graphical calculus induces a tensor functor
\(Z_{\catA}:\RT[\catA]\to\catA\), that is, graphical calculus for symmetric
categories is invariant by moves M\ref{M1}-M\ref{M2}.
\end{thm}

\subsection{Graphical calculus for ribbon graphs}
\label{sec:graph-calc-streb}

We want now to define a variant of graphical calculus suitable for
application to ribbon graphs. Such graphs arise both in Feynman
expansions of matrix integrals and in orbifold cellularizations of
moduli spaces of curves\footnote{These two subjects are indeed deeply
  related: see \cite{di-francesco-itzykson-zuber;kontsevich-model,
    kontsevich;intersection-theory;1992, witten;kontsevich-model}.}.
Any ribbon graph can be realized (in many different ways) as an
RT-graph. So, in order to define a graphical calculus on ribbon graphs
it suffices to define it on RT-graphs in a way that is
independent of the particular realization chosen.

\label{sec:ribbon-graphs}
\begin{dfn}
A \emph{ribbon graph} of type $(p,q)$ is a purely 1-dimensional
CW-complex \(\Gamma\) with 
\begin{enumerate}
\item $p+q$ endpoints divided into two ordered subsets
\(\In(\Gamma)\)
and \(\Out(\Gamma)\) with \(\card{\In(\Gamma)}=p\) and
\(\card{\Out(\Gamma)}=q\);
\item a cyclic order on half edges stemming from each vertex.
\end{enumerate}
The linear spans $\RG(p,q)$ of ribbon graphs of type $(p,q)$ define the
\(\Hom\)-spaces of the PROP $\RG$ of ribbon graphs.
\end{dfn}

\begin{figure}[htbp]
  \usegraphics{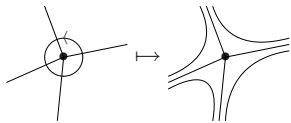}{%
    \begin{equation*}
      {\xy
        ,(0,0)*{\bullet},(0,0);(1,0.2)**\dir{-},(0,0);(-0.3,0.8)**\dir{-}%
        ,(0,0);(-0.9,-0.4)**\dir{-},(0,0);(-0.1,-1)**\dir{-}%
        ,(0,0.3);(0,0)%
        ,{\ellipse<>{-}},(0.02,0.3)*{\scriptscriptstyle{\langle}}
        \endxy}
      \mapsto
      {\xy
        ,(0,0)*{\bullet},(0,0);(1,0.2)**\dir{-},(0,0);(-0.3,0.8)**\dir{-}%
        ,(0,0);(-0.9,-0.4)**\dir{-},(0,0);(-0.1,-1)**\dir{-}%
        ,(1,0.1);(0,-1)**\crv{(0.05,-0.05)}%
        ,(-0.2,0.8);(1,0.3)**\crv{(0.05,0.05)}%
        ,(-0.9,-0.3);(-0.4,0.8)**\crv{(-0.05,0.05)}%
        ,(-0.9,-0.5);(-0.2,-1)**\crv{(-0.05,-0.05)}%
        \endxy}
    \end{equation*}
    }
  \caption{Fattening edges at a vertex with cyclic order.}
  \label{fig:fattening-edges} 
\end{figure}

Ribbon graphs arose in connection with a certain cellular
decomposition of the moduli space of smooth complex curves (see
\cite{harer;cohomology-of-moduli, mulase-penkava}). The connection is,
very roughly, the following: choose a ribbon graph $\Gamma$ of type
$(0,0)$; one can use the cyclic order to ``fatten'' edges
into thin ribbons\footnote{Hence the name ``ribbon graph''.} (see
\prettyref{fig:fattening-edges}) --- so we turn the graph into a
compact oriented surface with boundary \(S(\Gamma)\); this
construction may be refined to take into account a conformal structure
on $S(\Gamma)$. The boundary components of $S(\Gamma)$ are called
``holes'' of the ribbon graph \(\Gamma\); the set of holes of
\(\Gamma\) is denoted $\Holes{\Gamma}$.
The number of boundary components \(s\) and the genus \(g\) of the ribbon
graph \(\Gamma\) are defined to be those of the surface \(S(\Gamma)\).
Notice that it is meaningful to speak of the genus and number of boundary
components only for ribbon graphs of type (0,0).

There is a natural forgetful functor \(\RT\to\RG\) which forgets
orientations on edges and at each vertex \(v\) remembers only the cyclic
order
induced by the total order on \(\In(v)\) and \(\Out(v)\). 
\begin{equation*}
  \usegraphics{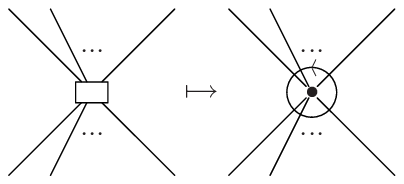}{%
    \xy*!LC\xybox{
      (0,1)*+[F]{\ };%
      (-1,0)**\dir{-},(-0.5,0)**\dir{-},%
      (0,0.5)*+{\ldots},(1,0)**\dir{-},%
      (-1,2)**\dir{-},(-0.5,2)**\dir{-},%
      (0,1.5)*+{\ldots},(1,2)**\dir{-},%
      }\endxy\mapsto
    \xy*!LC\xybox{
      (0,1)*{\bullet};%
      (-1,0)**\dir{-},(-0.5,0)**\dir{-},%
      (0,0.5)*+{\ldots},(1,0)**\dir{-},%
      (-1,2)**\dir{-},(-0.5,2)**\dir{-},%
      (0,1.5)*+{\ldots},(1,2)**\dir{-},%
      ,(0,1.3);(0,1)%
      ,{\ellipse<>{-}},(0.02,1.3)*{\scriptscriptstyle{\langle}}
      }\endxy
    }
\end{equation*}
\begin{lemma}
The PROP \(\RG\) is the quotient of \(\RT\) with respect to relations
generated by the following moves
\begin{enumerate}[M1)]\setcounter{enumi}{2}
\item \label{M3} reverse orientation on edges;
\item \label{M4} the first (resp. the last) edge in \(\In(v)\) becomes the
first (resp. the last) edge in \(\Out(v)\) and vice-versa
\begin{equation*}
  \usegraphics{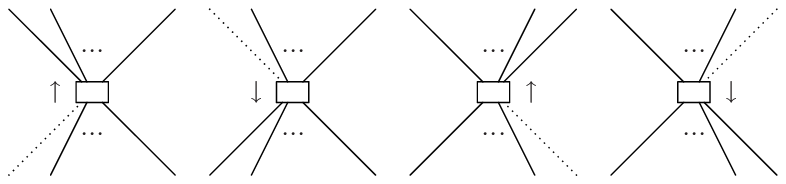}{%
    \xy*!LC\xybox{
      (0,1)*+[F]{\ };%
      (-1,0)**\dir{.},(-0.5,0)**\dir{-},%
      (0,0.5)*+{\ldots},(1,0)**\dir{-},%
      (-1,2)**\dir{-},(-0.5,2)**\dir{-},%
      (0,1.5)*+{\ldots},(1,2)**\dir{-},%
      (-.45,1)*{\uparrow}
      }\endxy\quad,\quad\xy*!LC\xybox{
      (0,1)*+[F]{\ };%
      (-1,0)**\dir{-},(-0.5,0)**\dir{-},%
      (0,0.5)*+{\ldots},(1,0)**\dir{-},%
      (-1,2)**\dir{.},(-0.5,2)**\dir{-},%
      (0,1.5)*+{\ldots},(1,2)**\dir{-},%
      (-.45,1)*{\downarrow}
      }\endxy\quad,\quad
    \xy*!LC\xybox{
      (0,1)*+[F]{\ };%
      (-1,0)**\dir{-},(0.5,0)**\dir{-},%
      (0,0.5)*+{\ldots},(1,0)**\dir{.},%
      (-1,2)**\dir{-},(0.5,2)**\dir{-},%
      (0,1.5)*+{\ldots},(1,2)**\dir{-},%
      (.45,1)*{\uparrow}
      }\endxy\quad,\quad\xy*!LC\xybox{
      (0,1)*+[F]{\ };%
      (-1,0)**\dir{-},(0.5,0)**\dir{-},%
      (0,0.5)*+{\ldots},(1,0)**\dir{-},%
      (-1,2)**\dir{-},(0.5,2)**\dir{-},%
      (0,1.5)*+{\ldots},(1,2)**\dir{.},%
      (.45,1)*{\downarrow}
      }\endxy
    }
\end{equation*}
\end{enumerate}
\end{lemma} 

\begin{lemma}
  \label{lemma:passing-up-and-down}
In any rigid symmetric tensor category \(\catA\) there are
natural isomorphisms
\begin{equation*}
\catA(X\otimes Y, Z) \longleftrightarrow \catA(X, Z\otimes\rdual{Y}),
\qquad
\catA(X\otimes Y, Z) \longleftrightarrow \catA(Y, \ldual{X}\otimes Z),
\end{equation*}
for all \(X, Y, Z \in \catA\).
\end{lemma}
For instance, the second bijection in the above lemma maps a morphism
\(f: X \otimes Y \to Z\) to the composition
\begin{equation*}
X \xrightarrow{\coev_Y \otimes \id_X} \ldual{Y} \otimes Y \otimes X %
\xrightarrow{\id_{\ldual{Y}} \otimes f} \rdual{Y} \otimes Z.%
\end{equation*}
In graphical notation, this reads:
\begin{equation*}
  \usegraphics{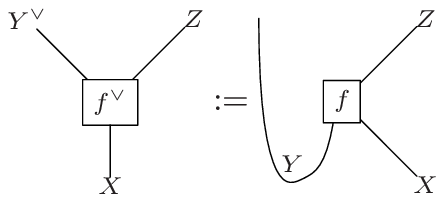}{\xy*!LC\xybox{%
    (0,1)*+[F]{\rdual{f}};%
    (0,0)*{X}**\dir{-},%
    (-1,2)*{\rdual{Y}}**\dir{-},(1,2)*{Z}**\dir{-},%
    }\endxy
  :=
  \xy*!LC\xybox{%
    (0,1)*+[F]{f}="x";%
    (1,0)*{X}**\dir{-},%
    \vloop~{(-1,0.5)}{"x"+DL+/d6mm/}{(-1,2)}{"x"+DL+/r1mm/}|{Y},%
    (1,2)*{Z}**\dir{-},%
    }\endxy}
\end{equation*}

Let \(V\) be a vector space over \(\fk\) equipped with a symmetric inner
product \(b: V\tp{2} \to \fk\). The category \(\langle V \rangle\) of
\prettyref{sec:vector-space-cat} is rigid symmetric monoidal; every
object in \(\langle V \rangle\) is self-dual. Therefore, bijections from
\prettyref{lemma:passing-up-and-down} translate into linear
isomorphisms:
\begin{gather}
\Hom(V\tp{p}\otimes V\tp{q}, V\tp{r}) \longleftrightarrow \Hom(V\tp{p},
V\tp{r}\otimes V\tp{q}),
\label{eq:right-rot}
\\
\Hom(V\tp{p}\otimes V\tp{q}, V\tp{r}) \longleftrightarrow \Hom(V\tp{q},
V\tp{p}\otimes V\tp{r}), 
\label{eq:left-rot}
\end{gather}
for all \(p,q,r \in \setZ\).

\begin{dfn}
A cyclic algebra structure on \((V, b)\) is a sequence
\(\{T_r\}_{r\in\setN}\) of cyclically invariant linear maps \(T_r :
V\tp{r} \to \fk\),:
\begin{equation*}
T_r (X_1 \otimes \dots \otimes X_{r-1} \otimes X_r) = 
T_r (X_r \otimes X_1 \otimes \dots \otimes X_{r-1}).
\end{equation*}
\end{dfn}
Fix a cyclic algebra structure \(\{T_r\}\) on \(V\).  Each map \(T_r\)
in turn defines linear maps \(T_{p,q}: V\tp{p} \to V\tp{q}\), for all
\(p\) and \(q\) such that \(p+q=r\), via the isomorphisms
\prettyref{eq:right-rot} and \prettyref{eq:left-rot}; cyclical
invariance guarantees that \(T_{p,q}\) does not depend on the
particular sequence of isomorphisms \prettyref{eq:right-rot} and
\prettyref{eq:left-rot}: any one yielding the right source and
target is good.

\begin{thm}\label{thm:gccyc}
Let \(V\) be a vector space. Then the following data are equivalent:
\begin{enumerate}
\item cyclic algebra structures on \(V\);
\item \(\RG\)-algebra structures on \(V\).
\end{enumerate}
\end{thm}
\begin{proof} Let \(Z:\RG\to\EndOp[V]\) be a PROP
action. Then 
\begin{align*}
b &:= Z\left(\xy*!LC\xybox{%
\vloop~{(0,0.5)}{(1,0.5)}{(0,-0.5)}{(1,-0.5)},(0,-0.7)*{1^{\textrm{in}}}%
,(1,-0.7)*{2^{\textrm{in}}}}
\endxy\right)
\\
T_r &:= Z\left(\xy*!LC\xybox{
(0,1)*{\bullet};%
(-1,0)*+{1^{\textrm{in}}}**\dir{-},(-0.5,0)*+{2^{\textrm{in}}}**\dir{-},%
(0,0.5)*+{\ldots},(1,0)*+{i^{\textrm{in}}}**\dir{-},%
(-1,2)*+{r^{\textrm{in}}}**\dir{-}%
,(-0.5,2)*+{\ \ \ (r-1)^{\textrm{in}}}**\dir{-},%
(0,1.5)*+{\ldots},(1,2)*+{(i+1)^{\textrm{in}}}**\dir{-},%
}\endxy\right) 
\end{align*}
defines a cyclic algebra structure on \(V\).
Vice-versa, let \((V,b,T_1,T_2,\dots)\) be a cyclic algebra. 
Pick any ribbon graph \(\Gamma\) and realize it as an RT-graph; call this
RT-graph \(\hat\Gamma\). Give \(\hat\Gamma\) the structure of a
\(\langle V\rangle\)-colored RT-graph by coloring all the edges of
\(\hat\Gamma\) with \(V\) and coloring any vertex \(v\) with
\(T_{\card{\In(v)},\card{\Out(v)}}\). Denote the \(\langle
V\rangle\)-colored RT-graph obtained this way by \(\hat\Gamma_V\). Two
realizations of \(\Gamma\)
as an RT-graph differ by a finite sequence of moves M3, M4. Since \(V\) is
self-dual, the graphical calculus \(Z_{\langle V\rangle}(\hat\Gamma_V)\)
is independent of the orientation on the edges, i.e., it  is invariant with
respect to the move M3. Moreover relations provided
by \prettyref{lemma:passing-up-and-down} give the invariance of
\(Z_{\langle V\rangle}(\hat\Gamma_V)\) with respect to moves of type M4:
\begin{equation*}
  \usegraphics{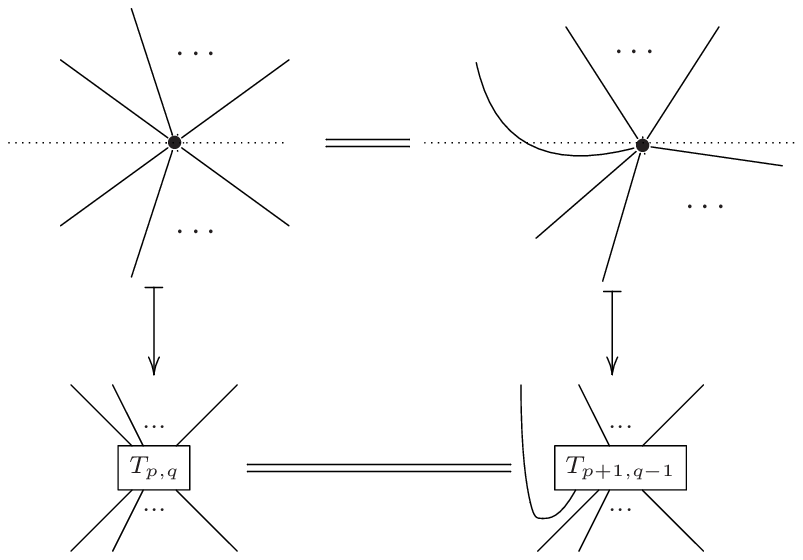}{%
    \xymatrix{
      \begin{xy}
        (0,0)*!LC\xybox{\rgvertex{10}%
          \loose2\missing3\loose4\loose5
          \loose7\loose8\missing9\loose{10}},
        (-0.5,0)*!LC{\ };(3,0)**\dir{.}
      \end{xy}
      \ar@{|->}[d] \ar@{=}[r]&      
      \begin{xy}
        (0,0)*!LC\xybox{\rgvertex{11}%
          \loose2\missing3\loose4,
          (-1,1);"CENTER"**\crv{"AUXii6"+/d6pt/},
          \loose7\loose8\missing{10}\loose{11}},
        (-1,0);(3.5,0)*!LC{\ }**\dir{.}
      \end{xy}
      \ar@{|->}[d]
      \\
      \xy*!LC\xybox{%
        (0,1)*+[F]{T_{p,q}};%
        (-1,0)**\dir{-},(-0.5,0)**\dir{-},(0,0.5)*{\ldots},(1,0)**\dir{-},%
        (-1,2)**\dir{-},(-0.5,2)**\dir{-},(0,1.5)*{\ldots},(1,2)**\dir{-},%
        }\endxy
      \ar@{=}[r]&
      \xy*!LC\xybox{%
        (0,1)*+[F]{T_{p+1,q-1}}="x";%
        (-1,0)**\dir{-},(-0.5,0)**\dir{-},(0,0.5)*{\ldots},(1,0)**\dir{-},%
        \vloop~{(-1.2,0.9)}{"x"+DL+/d11pt/}{(-1.2,2)}{"x"+DL+/r6pt/},%
        (-0.5,2)**\dir{-},(0,1.5)*{\ldots},(1,2)**\dir{-},%
        }\endxy
      }
    }
\end{equation*}
Therefore, \(Z(\Gamma):=Z_{\langle V\rangle}(\hat\Gamma_V)\) is well
defined and is a PROP action.
\end{proof}
\begin{rem}
  The tensors $T_r$ are \emph{generators} of the $\RG$-algebra
  structure on $V$; as such, they are independent and need not satisfy
  any further compatibility relation. For instance, cyclic algebras
  need not be associative.
\end{rem}

\subsection{Graphical calculus for ordinary graphs}
\label{sec:graph-calc-graphs}

It is now easy to adapt constructions above to ordinary graphs. 
Note that a graph can be obtained from a ribbon graph by forgetting
the cyclic order on the edges incident to any vertex. Two ribbon
graphs leading to the same graph differ just by a permutation of the
edges incident to some vertices, so, to define a graphical calculus
for graphs, it will suffice to have a graphical calculus for ribbon
graphs, invariant with respect to the action of symmetric groups.

\begin{dfn}
  \label{dfn:symmetric-algebra}
  A symmetric algebra structure on \((V, b)\) is a sequence
  \(\{S_r\}\) of linear maps \(S_r : V\tp{r} \to \fk\) such that
  \begin{equation*}
    S_r (X_1 \otimes \dots \otimes X_r) = S_r
    (X_{\sigma(1)} \otimes\dots \otimes X_{\sigma(r)}), 
    \quad \forall\sigma\in\Perm{r},
  \end{equation*}
  that is, maps $\{S_r\}$ are invariant with respect to the action of
  the symmetric group.
\end{dfn}
Like in a cyclic algebra, the tensors $S_r$ are independent one from
the other: they should be regarded as generators of a PROP-algebra
structure.
\begin{dfn}
  \label{dfn:symmetric-graph-category}
  Denote by \(\CG\) the PROP of (ordinary topological) graphs:
  \(\CG(p,q)\) is the linear span of graphs of type $(p,q)$. It is the
  quotient of \(\RG\) by the action of symmetric groups on half-edges
  stemming from any vertex.
\end{dfn}  

\begin{thm}\label{thm:gcsym}
Let \(V\) be a vector space. Then the following data are equivalent:
\begin{enumerate}
\item symmetric algebra structures on \(V\);
\item \(\CG\)-algebra structures on \(V\).
\end{enumerate}
\end{thm}
\begin{proof}
An action \(Z:\CG\to\EndOp[V]\) endows \(V\) with a symmetric 
algebra structure, as in the proof of
\prettyref{thm:gccyc}. Conversely, assume we are given a cyclic algebra
structure on \(V\). Since the
tensors $S_r$ are symmetric, they are, in particular, cyclically invariant:
so
\((V,b,S_1,S_2,\dots)\) is a cyclic algebra and an action
\(Z:\RG\to\EndOp[V]\) is defined. 
Since the tensors \(\{S_k\}\) are symmetric, this action
factors through an action \(Z:\CG\to\EndOp[V]\).
\end{proof}
\begin{rem}
\label{rem:many-sorted-graphs}
So far, we have considered cyclic (resp. symmetric) algebras with
only one $r$-ary operation for any $r\in\setN$. The graphical calculus
formalism immediately generalizes to cyclic (resp. symmetric)
algebras with a family of cyclic tensors \(\{T_{r,\alpha}\}_{\alpha\in
I_r}\)
(resp. a family of symmetric tensors \(\{S_{r,\alpha}\}_{\alpha\in
I_r}\));
in fact, we only ought to consider graphs whose $r$-valent vertices
are decorated with labels from $I_r$.

The \emph{dual graphs} of stable curves \cite{deligne-mumford} provide
an example --- they are ordinary graphs, with each vertex \(v\)
decorated by an integer \(g(v)\): it is the genus of an irreducible
component of the algebraic curve corresponding to the graph. Such
graphs were called \emph{modular graphs} in \cite{getzler-kapranov}.
One can apply methods of graphical calculus to this class of graphs;
we give a specimen in \prettyref{xmp:getzler-formula}.
\end{rem}

\subsection{A sample computation}
Let \((V,b,S_1,S_2,\dots)\) be a symmetric algebra. We want to compute
the operator
\begin{equation*}
Z(\Gamma):=Z\left(
\xy
\rghole{3}\loose{1}\loose{2}\loose{3},(0,2)*{1^{\text{in}}}%
,(1.73,-1)*{2^{\text{in}}},(-1.73,-1)*{3^{\text{in}}}
\endxy
\right):V^{\otimes 3}\to\fk
\end{equation*}
A realization of the graph \(\Gamma\) as an RT-diagram is obtained by
hanging it by the vertices:
\begin{equation*}
  \usegraphics{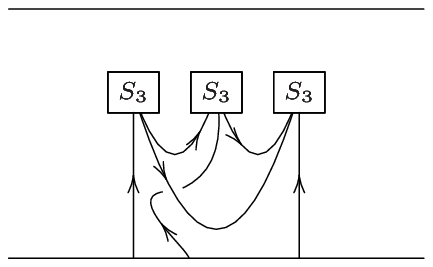}{\xy*!LC\xybox{
      (-1,1)="a",(0,1)="b",(1,1)="c",(-0.4,-0.15)="buno",(-0.65,-0.2)="bdue",%
      (-1,-1)="d",(-.33,-1)="e",(1,-1)="f",%
      "a"*+[F]{S_3};"d"**\dir{-}?(.5)*\dir{<},%
      "c"*+[F]{S_3};"f"**\dir{-}?(.5)*\dir{<},%
      "b"*+[F]{S_3};"buno"**\crv{(0.15,0.1)},%
      "bdue";"e"**\crv{(-1.1,-0.27)&(-0.33,-0.9)}?(.5)*\dir{<},%
      "a"*+[F]{S_3};"b"*+[F]{S_3}**\crv{(-.5,-.5)}?(.8)*\dir{>},%
      "a"*+[F]{S_3};"c"*+[F]{S_3}**\crv{(0,-2.3)}?(.2)*\dir{>},%
      "b"*+[F]{S_3};"c"*+[F]{S_3}**\crv{(0.5,-.5)}?(.3)*\dir{>},%
      (-2.5,2);(2.5,2)**\dir{-},(-2.5,-1);(2.5,-1)**\dir{-}
      }\endxy}
\end{equation*}
Splitting the RT-diagram above into elementary pieces and applying
Reshetikhin-Turaev rules we find
\begin{multline*}
  \usegraphics{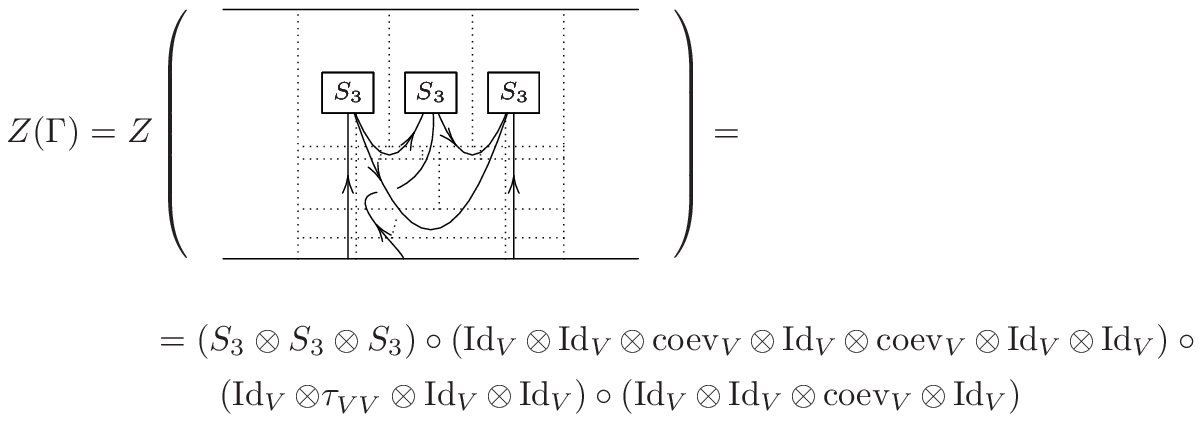}{%
    Z(\Gamma) = Z\left(\quad{\xy*!LC\xybox{
          (-1,1)="a",(0,1)="b",(1,1)="c",(-0.4,-0.15)="buno",(-0.65,-0.2)="bdue",%
          (-1,-1)="d",(-0.33,-1)="e",(1,-1)="f",%
          "a"*+[F]{S_3};"d"**\dir{-}?(.5)*\dir{<},%
          "c"*+[F]{S_3};"f"**\dir{-}?(.5)*\dir{<},%
          "b"*+[F]{S_3};"buno"**\crv{(0.15,0.1)},%
          "bdue";"e"**\crv{(-1.1,-0.27)&(-0.33,-0.9)}?(.5)*\dir{<},%
          "a"*+[F]{S_3};"b"*+[F]{S_3}**\crv{(-.5,-.5)}?(.8)*\dir{>},%
          "a"*+[F]{S_3};"c"*+[F]{S_3}**\crv{(0,-2.3)}?(.2)*\dir{>},%
          "b"*+[F]{S_3};"c"*+[F]{S_3}**\crv{(0.5,-.5)}?(.3)*\dir{>},%
          (-2.5,2);(2.5,2)**\dir{-},(-2.5,-1);(2.5,-1)**\dir{-},
          (-1.6,0.35);(1.6,0.35)**\dir{.},%
          (-1.6,0.2);(1.6,0.2)**\dir{.},%
          (-1.6,-0.4);(1.6,-0.4)**\dir{.},%
          (-1.6,-0.75);(1.6,-0.75)**\dir{.},%
          (-1.6,2);(-1.6,-1)**\dir{.},%
          (1.6,2);(1.6,-1)**\dir{.},%
          (-0.5,2);(-0.5,0.35)**\dir{.},%
          (0.5,2);(0.5,0.35)**\dir{.},%
          (-0.1,0.35);(-0.1,0.2)**\dir{.},%
          (0.1,0.35);(0.1,-0.4)**\dir{.},%
          (-0.9,0.8);(-0.9,-1)**\dir{.},%
          (0.9,0.8);(0.9,-1)**\dir{.},%
          (-0.6,0.35);(-0.6,0.2)**\dir{.},%
          (0.6,0.35);(0.6,0.2)**\dir{.},%
          (-0.4,-0.47);(-0.47,-0.75)**\dir{.}}\endxy}\quad\right)
    \\
    = (S_3 \otimes S_3 \otimes S_3) \circ (\Id_V \otimes \Id_V \otimes
    \coev_V \otimes \Id_V \otimes \coev_V \otimes \Id_V \otimes \Id_V)
    \circ
    \\
    (\Id_V \otimes \tau_{VV}^{} \otimes \Id_V \otimes \Id_V) \circ
(\Id_V
    \otimes \Id_V \otimes \coev_V \otimes \Id_V) }
\end{multline*}
If \(\{e_i\}\) is a basis of \(V\) and \(\{e^i\}\) denotes the dual
basis with respect to the pairing \(b\), then structure constants of
the symmetric algebra \((V,b,S_1,S_2,\dots)\) are
\begin{equation*}
g_{ij}=b(e_i,e_j),\qquad g^{ij}:=b(e^i,e^j),
\end{equation*}
\begin{equation*}
(S_k)_{i_1,\dots,i_s}^{j_{s+1},\dots,j_k}:=S_k(e_1\otimes\cdots\otimes
e_{i_s}\otimes e^{j_{s+1}}\otimes\cdots\otimes e^{j_k}).
\end{equation*}
With these notations
\begin{equation*}
\coev_V=\sum_{i,j}g^{ij}e_i\otimes e_j,
\end{equation*}
and the operator \(Z(\Gamma)\) acts on basis elements by
\begin{equation*}
Z(\Gamma)(e_\alpha\otimes e_\beta\otimes
e_\gamma)=\sum_{\delta,\epsilon,\zeta,\eta,\vartheta,\iota}g^{\delta\epsilon}
g^{\zeta\eta}g^{\vartheta\iota} (S_3)_{\alpha\delta\zeta}
(S_3)_{\eta\beta\vartheta} (S_3)_{\iota\epsilon\gamma}.
\end{equation*}

\subsection{On graphical notations used in physics literature}
All the graphs we have considered so far are grouped under the generic name
of ``Feynman diagrams''; the Casimir element and the tensors
decorating the vertices are called ``propagator'' and
``interactions'', respectively. Moreover, different types of lines
are used to denote different kinds of particles (see, for instance
\cite{itzykson-zuber;quantum-field-theory, QFS}); one can recover
these notations as follows.

A graph \(\Gamma\) of type \((r,0)\) gives a linear operator
\(Z(\Gamma):V\tp{r}\to\setC\).  In graphical notations, the value of this
operator on \(v_1\otimes v_2\otimes\cdots\otimes v_r\) can be represented
by the graph \(\Gamma\) with the \(i\)-th incoming edge decorated by the
vector \(v_i\): 
\begin{equation*} Z\left( \xy ,(0,0)*{\bullet}%
,(0,0);(1.8,0)**\dir{-},(0,0);(-0.7,-0.7)**\dir{-}%
,(0,0);(-0.7,0.7)**\dir{-},(1.8,0);(2.5,-0.7)**\dir{-}%
,(1.8,0);(2.5,0.7)**\dir{-}%
,(1.8,0)*{\bullet}
,(-0.9,-0.9)*{1^{\textrm{in}}},(-0.9,0.9)*{2^{\textrm{in}}}%
,(2.7,0.9)*{3^{\textrm{in}}},(2.7,-0.9)*{4^{\textrm{in}}} \endxy
\right)(v_1\otimes v_2\otimes v_3\otimes v_4)=:\xy ,(0,0)*{\bullet}%
,(0,0);(1.8,0)**\dir{-},(0,0);(-0.7,-0.7)**\dir{-}%
,(0,0);(-0.7,0.7)**\dir{-},(1.8,0);(2.5,-0.7)**\dir{-}%
,(1.8,0);(2.5,0.7)**\dir{-}%
,(1.8,0)*{\bullet}
,(-0.9,-0.9)*{v_1^{\textrm{in}}},(-0.9,0.9)*{v_2^{\textrm{in}}}%
,(2.7,0.9)*{v_3^{\textrm{in}}},(2.7,-0.9)*{v_4^{\textrm{in}}} 
\endxy
\end{equation*} 
The vectors \(v_i\) are called the ``incoming states''.  If a basis
\(\{e_i\}\) (not necessarily orthonormal) is given for \(V\), then an
incoming state \(e_i\) will be denoted simply by the index
\(i\).\footnote{This introduces a notations clash, since we use
  indices near edges to indicate the total order of $\In(\Gamma)$ and
  $\Out(\Gamma)$.}  If \(\{e^i\}\) is the dual basis of \(\{e_i\}\) with
respect to the inner product of \(V\), an incoming state \(e^i\) (or
an outgoing $e_i$) will be denoted by the index \(\rdual{i}\). It is
customary to write
\begin{equation*} 
\xy*!LC\xybox{
(0,1)*{\bullet};%
(-1,0)*+{i_1^{\textrm{in}}}**\dir{-},(-0.5,0)*+{i_2^{\textrm{in}}}**\dir{-},%
(0,0.5)*+{\ldots},(1,0)*+{i_r^{\textrm{in}}}**\dir{-},%
(-1,2)*+{j_1^{\textrm{out}}}**\dir{-},(-0.5,2)*+{\ \ %
j_2^{\textrm{out}}}**\dir{-},%
(0,1.5)*+{\ldots},(1,2)*+{j_s^{\textrm{out}}}**\dir{-},%
}\endxy:= \left(
\xy*!LC\xybox{ (0,1)*{\bullet};%
(-1,0)*+{i_1^{\textrm{in}}}**\dir{-},(-0.5,0)*+{i_2^{\textrm{in}}}**\dir{-},%
(0,0.5)*+{\ldots},(1,0)*+{i_r^{\textrm{in}}}**\dir{-},%
(-1,2)*+{{\rdual{j}_1}^{\textrm{in}}}**\dir{-}%
,(-0.5,2)*+{\ \ {\rdual{j}_2}^{\textrm{in}}}**\dir{-},%
(0,1.5)*+{\ldots},(1,2)*+{{\rdual{j}_s}^{\textrm{in}}}**\dir{-},%
}\endxy\right)\cdot e_{j_1}\otimes\cdots\otimes e_{j_s}. 
\end{equation*} 
For example, if \(T:V\tp{3}\to\setC\) is a cyclic or
symmetric tensor, we have: 
\begin{equation*} 
\xy ,(0,0)*{\bullet}%
,(0,0);(1.8,0)**\dir{-},(0,0);(-0.7,-0.7)**\dir{-}%
,(0,0);(-0.7,0.7)**\dir{-},(1.8,0);(2.5,-0.7)**\dir{-}%
,(1.8,0);(2.5,0.7)**\dir{-}%
,(1.8,0)*{\bullet}
,(-0.9,-0.9)*{i^{\textrm{in}}},(-0.9,0.9)*{j^{\textrm{in}}}%
,(2.7,0.9)*{k^{\textrm{in}}},(2.7,-0.9)*{l^{\textrm{in}}}
\endxy=\sum_mT_{ijm}T^m_{kl}=\sum_{m,n}T_{ijm}g^{mn}T_{nkl}
\end{equation*} 
When the space \(V\) of ``physical states'' is the direct
sum of two subspaces \(V_1\), \(V_2\), vectors are usually depicted by
different types of lines according to the subspace they lie in. For
example, it is customary in physical literature to depict \emph{fermions}
by a straight line and \emph{bosons} by a wavy line\footnote{ The
graphical convention above is often refined depicting spin 0 bosons by
dashed straight lines and reserving wavy lines for spin 1 bosons (see, for
example \cite{cheng-oneill}).}, so that one encounters diagrams like the
following, which depicts a photon exchange between two electrons:
\begin{equation*} 
\xy ,(0,0)*{\bullet}%
,(0,0);(1.8,0)**\dir{~},(0,0);(-0.7,-0.7)**\dir{-}%
,(0,0);(-0.7,0.7)**\dir{-},(1.8,0);(2.5,-0.7)**\dir{-}%
,(1.8,0);(2.5,0.7)**\dir{-},(1.8,0)*{\bullet}%
,(-0.9,-0.9)*{{(e^-)}^{\textrm{in}}},(-0.9,0.9)*{{(e^-)}^{\textrm{out}}}%
,(2.7,0.9)*{{(e^-)}^{\textrm{out}}},(2.7,-0.9)*{{(e^-)}^{\textrm{in}}}
\endxy 
\end{equation*}


\begin{rem}[Fields of algebras and the WDVV equation]
Let \(\phi\) be an analytic function defined on a neighborhood \(U\) of
\(0\in V\). Then, for any \(x\in U\), the derivatives
\(D^n\phi\vert_x\) are a family
of symmetric tensors on the tangent spaces \(T_xV\). These tensors define
a \emph{field} of symmetric algebras \(A_x\) on \(U\): \(A_x\) is  
the symmetric algebra
\((T_xV,\inner{-}{-},D\phi\vert_x,D^2\phi\vert_x,\dots)\) --- the
inner product is induced by the canonical identification of \(V\) with
its tangent spaces. The function \(\phi\) is called a \emph{potential}
for the field of algebras \(A_x\). Equivalently, the field of algebras
\(A_x\) is the datum of a symmetric algebra structure on the
\(C^\infty(U)\)-module \({\mathcal X}(U)\) of smooth vector fields on
\(U\). 

Denote by \(Z_{\phi,x}\) the graphical calculus for \(A_x\). Note that for
every graph \(\Gamma\) the map
\(\phi\mapsto Z_{\phi,x}(\Gamma)
\) is a differential operator --- we denote it by the symbol \(D_\Gamma\),
i.e., 
\begin{equation}
D_\Gamma(\phi):=Z_{\phi,x}(\Gamma),\quad\forall\phi\text{ analytic in }U.
\end{equation}
Abusing notation, we will occasionally write \(\Gamma(\phi)\) to mean
\(D_\Gamma(\phi)\); it comes handy when one is dealing with a field of
algebras enjoying some special property that can be described in
diagrammatic form. 

For example, symmetric associative algebras are described by the
associativity equation
\begin{equation*}
\xy
,(0,0)*{\bullet}%
,(0,0);(1.8,0)**\dir{-},(0,0);(-0.7,-0.7)**\dir{-}%
,(0,0);(-0.7,0.7)**\dir{-},(1.8,0);(2.5,-0.7)**\dir{-}%
,(1.8,0);(2.5,0.7)**\dir{-}%
,(1.8,0)*{\bullet}
,(-0.9,-0.9)*{x^{\textrm{in}}}%
,(-0.9,0.9)*{y^{\textrm{in}}}%
,(2.7,0.9)*{z^{\textrm{in}}}%
,(2.7,-0.9)*{w^{\textrm{in}}}
\endxy
=
\xy
,(0,-0.9)*{\bullet}%
,(0,-0.9);(0,0.9)**\dir{-},(0,-0.9);(-0.7,-1.6)**\dir{-}%
,(0,-0.9);(0.7,-1.6)**\dir{-},(0,0.9);(-0.7,1.6)**\dir{-}%
,(0,0.9);(0.7,1.6)**\dir{-}%
,(0,0.9)*{\bullet}
,(-0.9,-1.8)*{x^{\textrm{in}}}%
,(-0.9,1.8)*{y^{\textrm{in}}}%
,(0.9,1.8)*{z^{\textrm{in}}},(0.9,-1.8)*{w^{\textrm{in}}}
\endxy
\quad\forall\ x,y,z,w\in A.
\end{equation*}
So, the condition \(\phi\) must satisfy to define a field of symmetric
associative algebras is the Witten-Djikgraaf-Verlinde-Verlinde equation
(WDVV for short)
\begin{equation*}
\left(
\xy
,(0,0)*{\bullet}%
,(0,0);(1.8,0)**\dir{-},(0,0);(-0.7,-0.7)**\dir{-}%
,(0,0);(-0.7,0.7)**\dir{-},(1.8,0);(2.5,-0.7)**\dir{-}%
,(1.8,0);(2.5,0.7)**\dir{-}%
,(1.8,0)*{\bullet}
,(-0.9,-0.9)*{X^{\textrm{in}}}%
,(-0.9,0.9)*{Y^{\textrm{in}}}%
,(2.7,0.9)*{Z^{\textrm{in}}},(2.7,-0.9)*{W^{\textrm{in}}}
\endxy
\right)(\phi)=  
\left(\xy
,(0,-0.9)*{\bullet}%
,(0,-0.9);(0,0.9)**\dir{-},(0,-0.9);(-0.7,-1.6)**\dir{-}%
,(0,-0.9);(0.7,-1.6)**\dir{-},(0,0.9);(-0.7,1.6)**\dir{-}%
,(0,0.9);(0.7,1.6)**\dir{-}%
,(0,0.9)*{\bullet}
,(-0.9,-1.8)*{X^{\textrm{in}}}%
,(-0.9,1.8)*{Y^{\textrm{in}}}%
,(0.9,1.8)*{Z^{\textrm{in}}},(0.9,-1.8)*{W^{\textrm{in}}}
\endxy\right)(\phi)
,\quad\forall\ X,Y,Z,W\in {\mathcal X}(U)
\end{equation*} 
In terms of the flat vector fields
\(\{\partial_i\}\) on \(U\) corresponding to a basis \(\{e_i\}\) of
\(V\), the WDVV equation reads
\begin{equation}
\sum_{m,n}(\partial_i\partial_j\partial_m\phi)g^{mn}
(\partial_n\partial_k\partial_l\phi)=
\sum_{m,n}(\partial_i\partial_l\partial_m\phi)g^{mn}
(\partial_n\partial_j\partial_k\phi)
\end{equation}
which is the form one usually finds in the literature (e.g. 
\cite{kontsevich-manin}).
\end{rem}


\everyxy={0,<2em,0em>:,(0,0.5),} 

\section{Gaussian integrals and Feynman diagrams}
\label{sec:feynman-diagrams}

In this section we show how a Gaussian integral can be expanded into a sum
of Feynman diagrams, to be evaluated according to the rules of
graphical calculus. Depending on the nature of the integral, this formula
will hold as a strict equality or in the sense of asymptotic
expansions. In particle physics, Gaussian integrals and their Feynman
diagrams expansions are used to describe bosonic statistics. 
\subsection{Gaussian measures and the
  Wick's lemma} 
Let \(V\) be a finite dimensional Euclidean space,
with inner product \(\inner{-}{-}\). If \(\{e_i\}\) is a basis of
\(V\), we denote the \emph{coordinate} maps relative to this basis as
\(e^i:V\to\setR\), and write \(v^i\) for the pairing \(\pairing{e^i}{v}\).
The matrix associated to \(\inner{-}{-}\) with respect to the basis
\(\{e_i\}\) is given by
\begin{equation*}
  g_{ij} := \inner{e_i}{e_j}.
\end{equation*}
As customary, we set \(g^{ij} := (g^{-1})_{ij} = \inner{e^i}{e^j}\).

Let now \(\ud v\) be a (non trivial) translation invariant measure on
\(V\). The function \(\E^{-\onehalf\inner{v}{v}}\) is positive and
integrable with respect to \(\ud v\).
\begin{dfn}
  The probability measure on \(V\) defined by
  \begin{equation*}
    \ud\mu(v) := \frac{1}{A} \E^{-\onehalf\inner{v}{v}} \ud v, 
    \qquad
    A = {\int_V \E^{-\onehalf\inner{v}{v}} \ud v},
  \end{equation*}
  is called the \emph{Gaussian measure} on \(V\).
\end{dfn}
Since a non-trivial translation invariant measure on \(V\) is unique up
to a scalar factor, \(\ud\mu\) is actually independent of the chosen
\(\ud v\).

The symbol \(\avg{f}\) denotes the \emph{average} of a function \(f\)
with respect to the Gaussian measure, i.e.,
\begin{equation*}
  \avg{f} := \int_V f(v)\ud\mu(v).
\end{equation*}

\begin{lemma}[Wick]
  Polynomial functions of the coordinates \(v^i\) are integrable with
  respect to \(\ud\mu\) and:
  \begin{align}
    \tag{W1}
    \avg{v^{i_1}v^{i_2}\cdots
      v^{i_{2n+1}}} &= 0,
    \\
    \tag{W2}
    \avg{v^{i}v^{j}} &= g^{ij},
    \\
    \tag{W3}
    \avg{v^{i_1}v^{i_2}\cdots
      v^{i_{2n}}} &= \sum\limits_{s\in P}
    g^{i_{s_1} i_{s_2}} g^{i_{s_3} i_{s_4}} \cdots
    g^{i_{s_{2n-1}} i_{s_{2n}}}, 
  \end{align}
  where the sum ranges over all distinct pairings of the set of
  indices \(\{i_1,\dots,i_{2n}\}\), i.e., over the set of all
  partitions \(\{\{i_{s_1},i_{s_2}\},\{i_{s_3},i_{s_4}\},\dots\}\) of
  \(\{i_1,i_2,\dots,i_{2n}\}\) into 2-element subsets.
\end{lemma}
For a proof of Wick's lemma see, for instance,
\cite{bessis-itzykson-zuber;graphical-enumeration}.

The inner product \(\inner{-}{-}\) extends uniquely to a Hermitian
product on the complex vector space \(\cplx{V} := V \otimes \setC\).
Identify \(V\) with the subspace \(V \otimes \{1\}\) of real vectors in
\(\cplx{V}\); \(\{ e_i \}\) is a real basis for the complex vector
space \(\cplx{V}\). Extend \(\avg{-}\) to tensor powers of real
vectors by
\begin{equation*}
  \label{eq:avg-x}
  \avg{v\tp{k}} := \sum \avg{v^{i_1} \cdot \cdots \cdot v^{i_k}}
  e_{i_1} \otimes \cdots \otimes e_{i_k},
\end{equation*}
so Wick's lemma can be recast this way:
\begin{align}
  \label{eq:W1}\tag{W1'}
  \avg{v\tp{ 2n+1}} &= 0,
  \\
  \label{eq:W2}\tag{W2'}
  \avg{v\tp{ 2}} &= \sum_{i,j} g^{ij} e_i
  \otimes e_j,
  \\
  \label{eq:W3}\tag{W3'}
  \avg{v\tp{2n}} &= \sum\limits_{i_1, \dots,
    i_{2n}} \sum\limits_{s\in P} g^{i_{s_1} i_{s_2}} \cdots
  g^{i_{s_{2n-1}} i_{s_{2n}}} e_{i_1} \otimes e_{i_2} \otimes \cdots
  \otimes e_{i_{2n-1}} \otimes e_{i_{2n}},
\end{align}
where the last sum ranges over all distinct pairings of indices in the
set \(\{i_1, \dots, i_{2n}\}\).

The right-hand side of (\ref{eq:W2}) is the Casimir element
$\gamma_{\cplx{V}}$ of \(\{\cplx{V}, \inner{-}{-}\}\); in
Reshetikhin-Turaev's graphical notation, we can rewrite
\prettyref{eq:W2} as
\begin{equation*}
  \avg{v\otimes v}= {\xy\vloop-\endxy}\ .
\end{equation*}
The graphical notation becomes particularly suggestive 
(and useful) when applied to \prettyref{eq:W3}:
\begin{equation}
  \label{eq:avg-to-casimir}
  \begin{split}
    \avg{v\tp{ 4}} &=
    \usegraphics[-2em]{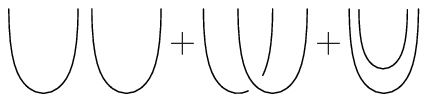}{%
      {\xy\vloop-,(1.2,0.5),\vloop-\endxy} +
      {\xy\vloop-|(.7)\knothole,(0.5,0.5),\vloop-\endxy} + 
      {\xy\vloop-,(0.7,0):(0.2,0.7),\vloop-\endxy}
    }\quad, \\
    \ldots& \\
    \avg{v\tp{ 2n}} &=
    \usegraphics[-2em]{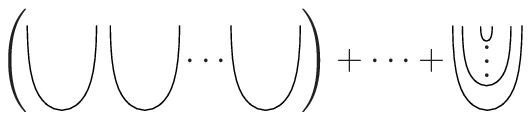}{%
      \left({\xy\vloop-,(1.2,0.5),\vloop-,(1.2,0)\endxy} \cdots
        {\xy\vloop-\endxy}\right) +
      \dots +
      {\xy\vloop-,(0.7,0):(0.2,0.7),%
        \vloop-,(0.7,0.2)*\txt{\vdots},(0.24,0):(2.4,2.9)\vloop-\endxy}
      }\quad,
  \end{split}
\end{equation}
the last sum ranging over all possible configurations of \(n\) Casimir
elements and the braiding being the trivial one: \(x\otimes y\mapsto
y\otimes x\).  

In addition, assume we have, for any \(k\), a cyclically invariant
\(k\)-tensor
\begin{equation*}
  T_{k}:\cplx{V}{}\tp{k} \to \setC;
\end{equation*}
which has the graphical representation
\begin{equation*}
  \underbrace{\xy(0.3333,0):
    (0,0)="a",
    (2,0)="b",
    (4,0)="c",
    (6,0)="d",
    (3,4)="v",
    (3,5)="w",
    "v";"a"**\dir{-},
    "v";"b"**\dir{-},
    "v";"d"**\dir{-},
    "c"*\txt{\ldots},
    "w"*\txt{\(T_{k}\)},
    "v"*{\bullet}
    \endxy}_{k}
\end{equation*}

The data \((V, \inner{-}{-}, T_1, T_2, \dots)\) define a cyclic
algebra structure, so we can use Reshetikhin-Turaev's graphical
calculus \(\Gamma\mapsto Z(\Gamma)\) to compute averages
\begin{equation*}
  \bigl\langle T_1(v)^{l_1} T_2(v\tp{ 2})^{l_2} \cdots
  T_{k}(v\tp{ k})^{l_k} \bigr\rangle.
\end{equation*}
\begin{lemma}\label{lemma:feynman-avg-gc}
  Any average \(\bigl\langle T_1(v)^{l_1} T_2(v\tp{ 2})^{l_2}
  \cdots T_{k}(v\tp{ k})^{l_k} \bigr\rangle\) is a linear
  combination \(\sum \alpha_\Gamma Z(\Gamma)\) where \(\Gamma\) runs
  in the set of ribbon graphs of type (0,0) with \(l_i\) vertices of
valence \(i\), for \(i=1,\dots,k\).
\end{lemma}
\begin{proof}
  By linearity of the \(T_k\)'s and \eqref{eq:avg-x},
  \begin{equation*}
    \avg{T_1(v)^{l_1} \cdots T_{k}(v\tp{ k})^{l_k}} = 
    \bigl(T_1\tp{l_1} \otimes \cdots \otimes T_{k}\tp{l_k} \bigr)
    \avg{v\tp{\sum_i il_i}}.
  \end{equation*}
  If \(\sum il_i\) is odd, \(\avg{v\tp{\sum il_i}}\) is zero, and the set
  of graphs considered in the statement is empty. If \(\sum il_i\) is
  even, according to graphical calculus rules, $\bigotimes_{i=1}^k
  T_i\tp{l_i}$ corresponds to the juxtaposition of $l_1$ univalent
  vertices, $l_2$ bivalent vertices, etc., up to $l_k$ vertices of
  valence $k$. In this case, \prettyref{eq:avg-to-casimir} translates
  $\avg{v\tp{\sum il_i}}$ into edges connecting these vertices in all
  possible ways.
\end{proof}

\begin{xmp}
  For example,
  \begin{align*}
    \avg{\left(T_2(v\tp{ 2})\right)^2} &= 
    \usegraphics[-3em]{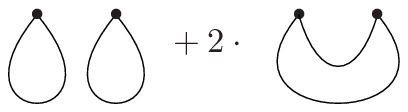}{%
      {\xy*!UC\xybox{
          /r 1cm/:
          ,(1,1.2);(1,1.2)**\crv{(0,0)&(2,0)},(1,1.2)*{\bullet}
          ,(1.8,1.2);(1.8,1.2)**\crv{(0.8,0)&(2.8,0)},(1.8,1.2)*{\bullet}
          }\endxy}
      + 2 \cdot 
      {\xy*!UC\xybox{
          /r 1cm/:
          ,(1,1.2);(1.8,1.2)**\crv{(0,0)&(2.8,0)},(1,1.2)*{\bullet}
          ,(1,1.2);(1.8,1.2)**\crv{(1.2,0.5)&(1.6,0.5)},(1.8,1.2)*{\bullet}
          }\endxy}
      }
    \\
    \avg{T_4(v\tp{ 4})} &=
    \usegraphics[-3em]{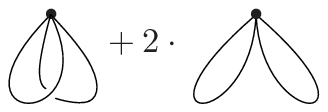}{%
      {\xy*!UC\xybox{
          /r 1cm/:
          ,(1,1.2);(1,1.2)**\crv{(-0.2,0)&(1.6,0)},(1,1.2)*{\bullet}
          ,(1,1.2);(1.05,.35)**\crv{(1,1.2)&(1.5,0.7)&(1.5,0.2)&(1.05,.35)}
          ,(1,1.2);(0.95,.45)**\crv{(1,1.2)&(0.8,0.5)&(0.95,.45)}
          }\endxy}
      + 2\cdot
      {\xy*!UC\xybox{
          /r 1cm/:
          ,(1,1.2);(1,1.2)**\crv{(-.4,0)&(.95,0)}
          ,(1,1.2);(1,1.2)**\crv{(1.05,0)&(2.4,0)}
          ,(1,1.2)*{\bullet}
          }\endxy}
      }
  \end{align*}
\end{xmp}

\begin{lemma}\label{lemma:feynman-avg-coeff}
  The coefficient \(\alpha_\Gamma\) appearing in
  \prettyref{lemma:feynman-avg-gc} is an integer given by:
  \begin{equation*}
    \alpha_\Gamma =
    \frac{1^{l_1} l_1! 2^{l_2} l_2! \cdots k^{l_k} l_k!}
    {\card{\Aut\Gamma}}.
  \end{equation*}
\end{lemma}
\begin{proof}
  Let \(X\) be the set of all ribbon graphs obtained by:
  \begin{inparaenum}
  \item juxtaposing \(l_1\) vertices of valence \(1\), \(l_2\)
    vertices of valence \(2\), etc., up to \(l_k\) vertices of valence
    \(k\), and,
  \item connecting them in all possible ways by means of arcs.
  \end{inparaenum}
  The constant \(\alpha_\Gamma\) counts the number of occurrences of graphs
  isomorphic to \(\Gamma\) in the set \(X\).
  
  The semi-direct product \(K = \prod_{i=1}^k(\Perm{l_i} \rtimes
(\setZ/i\setZ)^{l_i})\)
  acts on \(X\) as follows: the image of a graph
  \(\Phi\) is obtained by permuting vertices of the same valence and
  rotating edges incident to each vertex. Since this action is
  transitive on isomorphism classes,
  \begin{equation*}
    \alpha_\Gamma = \frac {\card{K}} {\card{\Stab_K(\Gamma)}}
    = \frac {\card{K}} {\card{\Aut\Gamma}}
    = \frac{1^{l_1} l_1! 2^{l_2} l_2! \cdots
      k^{l_k} l_k!}{\card{\Aut\Gamma}},
  \end{equation*}
  where $\Stab_K(\Gamma)$ is the stabilizer of $\Gamma$ under the
  action of $K$.
\end{proof}

\begin{thm}[Feynman-Reshetikhin-Turaev]\label{thm:FRT}
  Let \(Z_{x_*}\) be the graphical calculus for the cyclic algebra
  \(A(x_*) := (V_{\setC}, \inner{-}{-}, x_1T_1, x_2T_2, \dots)\), where
  $x_1,x_2,\dots$ are complex variables. Then, the following
  asymptotic expansion holds:
  \begin{equation}
    \label{eq:FRT1}
    Z(x_*):=\int_V \exp \biggl\{ \sum_{k=1}^\infty x_k
    \frac{T_k(v\tp{k})}{k}
    \biggr\} \ud\mu(v) 
    = \sum_{\Gamma\in\R(0,0)} \frac {Z_{x_*}(\Gamma)} {\card{\Aut\Gamma}},
  \end{equation}
  where the sum on the right ranges over the set $\R(0,0)$ of all
  isomorphism classes of (possibly disconnected) ribbon graphs of type
  (0,0). The formal series $Z(x_*)$ is called the \emph{partition
    function} of the algebra $A(x_*)$.
\end{thm}
\begin{proof}
  Expand in Taylor series the left-hand side:
  \begin{align*}
    \int_V \exp \biggl\{\sum_{k=1}^\infty &\frac{x_k}{k}T_k(v\tp{k})\biggr\}
    \ud\mu(v)=
    \\
    &= \sum_{n=0}^\infty \sum_{\nu_1, \dots, \nu_n} \frac{x_{\nu_1} \dots
      x_{\nu_n}} {n!\nu_1 \cdots \nu_n} \avg{T_{\nu_1}(v\tp{ \nu_1}) \cdots
      T_{\nu_n}(v\tp{\nu_n})}
    \\
    &= \sum_{k=0}^\infty \sum_{l_1, \dots, l_k} \frac{x_1^{l_1} \cdots x_k^{l_k}}
    {1^{l_1} l_1! 2^{l_2} l_2! \cdots k^{l_k} l_k!}  \avg{T_{1}(v)^{l_1}
      \cdots T_{k}(v\tp{k})^{l_k}}
    \\
    &= \sum_{k=0}^\infty \sum_{l_1, \dots, l_k} \frac
    {\avg{\bigl(x_1T_{1}(v)\bigr)^{l_1} \cdots
        \bigl(x_kT_{k}(v\tp{k})\bigr)^{l_k}}} {1^{l_1} l_1! 2^{l_2}
      l_2! \cdots k^{l_k} l_k!}, 
    \\
    &= \sum_{\Gamma \in \R(0,0)} \frac{Z_{x_*} (\Gamma)}
{\card{\Aut\Gamma}},
  \end{align*}
  by lemmas~\ref{lemma:feynman-avg-gc} and~\ref{lemma:feynman-avg-coeff}.
\end{proof}
A similar argument yields:
\begin{equation}
  \label{eq:FRT2}
  \int_V \frac{T_1(v)^{l_1} T_2(v\tp{2})^{l_2} \cdots
    T_k(v\tp{k})^{l_k}} {1^{l_1} l_1! 2^{l_2} l_2! \cdots k^{l_k}
    l_k!} \exp \left\{ \sum_{k=1}^\infty x_k\frac{T_k(v\tp{k})}{k}
    \right\} \ud\mu(v)
    = \sum_{\Gamma} \frac{Z_{x_*}(\Gamma)}{\card{\Aut\Gamma}}
\end{equation}
where the sum in the right-hand side ranges over all ribbon graphs
having $l_i$ ``special'' \(i\)-valent vertices (for $i=1\dots,k$), and
$\Aut(\Gamma)$ is the group of automorphisms that map the set of special
vertices to itself. In \prettyref{eq:FRT2}, the graphical calculus
$Z_{x_*}$ interprets each $i$-valent special vertex as the operator
$T_i$ and each ordinary $i$-valent vertex as the operator $x_iT_i$.
This is the same as considering graphs with two sorts of vertices (see
\prettyref{rem:many-sorted-graphs}), one decorated by operators $T_i$
(call them ``special vertices''), and the other decorated by $x_iT_i$
(call them ``ordinary'').


\prettyref{thm:FRT} can be straightforwardly adapted to a symmetric
algebra \((V, \inner{-}{-}, S_1,\break S_2, \dots)\).
\begin{thm}[Feynman-Reshetikhin-Turaev]\label{thm:FRTsym}
  Let \(Z_{x_*}\) be the graphical calculus for the symmetric algebra
  \((V_{\setC}, \inner{-}{-}, x_1S_1, x_2S_2, \dots)\), where
$x_1,x_2,\dots$ are
  complex variables. The following asymptotic expansion holds:
  \begin{equation*}
    Z(x_*):=\int_V \exp \biggl\{ \sum_{k=1}^\infty x_k
    \frac{S_k(v\tp{k})}{k!}
    \biggr\} \ud\mu(v)
    = \sum_{\Gamma\in\G(0,0)} \frac {Z_{x_*}(\Gamma)} {\card{\Aut\Gamma}},
  \end{equation*}
  where the sum on the right ranges over all isomorphism classes of
(possibly disconnected) ordinary graphs of type (0,0).
\end{thm}


\begin{rem}[Generalized Gaussian integrals]
  Denote by $Q(v)$ the Gaussian weight \(\E^{-1/2(v,v)}\), and let
  \(g^\sharp:V\to V^\lor\) be the isomorphism induced by the non-degenerate
  pairing \(g:V\otimes V\to\fk\). Moreover, let \(m\) be the usual
  multiplication on the space \(K := \fk[[V^\lor]]\) of formal power
  series on \(V\).  In \cite{oeckl;braided-qft}, Robert Oeckl shows
  that the graphical version of Wick's lemma is a formal consequence
  of the following properties:
\begin{enumerate}[(G1)]
\item\label{G1} a ``braided Leibniz rule'' for derivations:
\[\partial_w\circ m=m\bigl((\partial_w\otimes
\Id)+\tau^{-1}_{K, K}\circ(\partial_w\otimes\Id) \circ \tau^{}_{K,
  K}\bigr),\] for any \(w\in V\) and any \(\phi, \psi\in K\);
\item\label{G2} \(\partial_w Q = -g^\sharp(w)\cdot Q\), for all \(w\in V\);
\item\label{G3} \(g^\sharp\) is an isomorphism and
\(\ev_V\circ(\Id_V\otimes g^\#)=\ev_V\circ(\Id_V\otimes
g^\#)\circ\tau_{V,V}\);
\item\label{G4} 
\(\displaystyle{\int\partial_w(\phi\cdot Q)\ud
v=0}\), for
  any \(w\in V\) and any polynomial \(\phi\in K\);
\end{enumerate}
Equations (G\ref{G1}-G\ref{G4}) can be used to \emph{define} Gaussian
integrals in the context of \emph{arbirtary} braided monoidal categories.
This has been done by R.Oeckl, with the development of ``Braided
QFT'' \cite{oeckl;braided-qft}; since equation
\prettyref{eq:avg-to-casimir} is a formal consequence of   
(G\ref{G1}-G\ref{G4}), the whole machinery of Feynman diagrams expansions
will be available for these generalized Gaussian integrals, too.
A remarkable by-product of this general
theory is that the Berezin super-integrals of fermionic statistics are
simply obtained as ``braided
Gaussian integrals'' for a vector space \(V\) endowed with the non-trivial
braiding \(x\otimes y\mapsto -y\otimes x\) (the pairing \(g\) must,
consequently, be antisymmetric). Therefore, in the particular case
of the symmetric category of super vector spaces (with the usual graded
symmetric braiding), ``braided Gaussian integrals'' provide an unified
language for both statistics encountered in standard quantum field theory:
bosons correspond to even vectors and fermions correspond to odd vectors
\cite[Sections~3.3 and~3.4]{oeckl;spin-and-statistics}.
\end{rem}
%
%

\begin{xmp} 
Let \(\phi\) be an analytic function defined in a neighborhood of
\(0\in V\). The Taylor expansion formula  
\begin{equation*}
\phi(v)=\phi(0)+\sum_{n=1}^{\infty}\frac{D^n\phi\vert_0(v^{\otimes
n})}{n!},
\end{equation*}
together with \prettyref{thm:FRTsym}  gives
\begin{equation}
  \int_V \E^{\phi(v)} \ud
  \mu(v) = \E^{\phi(0)} \cdot \left(\sum_\Gamma
    \left.\frac{D_\Gamma(\phi)}{\card{\Aut\Gamma}}\right\vert_0\right).
\end{equation}
\end{xmp}

\begin{xmp}
  \label{xmp:getzler-formula}
  We want to show how this machinery can be used to derive a
  graph expansion formula cited in the introduction to
  \cite{getzler-kapranov}. If we are given several symmetric tensors
  \(\{S_{k,\alpha_k}\}_{\alpha_k\in I_k}\) for any index \(k\), the formula in
  \prettyref{thm:FRTsym} generalizes to
  \begin{equation}
    \label{eq:FRTeqter}
    \int_V \exp \biggl\{ \sum_{k=1}^\infty\sum_{\alpha_k\in I_k}
x_k{\nobreakspace}y_{\alpha_k}
    \frac{S_{k,\alpha_k}(v\tp{k})}{k!}   
    \biggr\} \ud\mu(v)
    = \sum_{\Gamma \in \G^I(0,0)} \frac {Z_{x_*,y_*}(\Gamma)}
{\card{\Aut\Gamma}}.
\end{equation}
where \(\G^I(0,0)\) denotes the set of isomorphism classes of (possibly
disconnected) graphs whose $k$-valent vertices are decorated by
elements of the set $I_k$. 

In particular, if we set, for any $k$,
\begin{equation*}
I_k=\mathbb{N}, \quad
x_k=1, \quad
y_g=\hbar^g, \quad
S_{k,g}=0\quad \text{ if } \quad 3g-3+k\leq0,
\end{equation*}
and rescale the inner product on $V$ by a factor $1/\hbar$ (i.e., we take
$(1/\hbar)\inner{-}{-}$ instead), then:
\begin{equation}
    \int_V \exp \biggl\{\frac{1}{\hbar} \sum_{k=1}^\infty\sum_{g=0}^\infty
\hbar^g
\frac{S_{k,g}(v\tp{k})}{k!}
    \biggr\} \ud\mu_\hbar(v)
    = \sum_{g=0}^\infty\hbar^g\sum_{\Gamma \in \textrm{M}(g,0)}
\hbar^{-b_0(\Gamma)}\frac
{Z(\Gamma)}
{\card{\Aut\Gamma}}.   
\end{equation}
where $\textrm{M}(g,0)$ denotes the set of isomorphism classes of
(possibly non connected) genus $g$ modular graphs with no legs. The
\emph{genus} of the modular graph $\Gamma$ is the integer $g(\Gamma):=\sum_v
g(v)+\dim H^1(\Gamma)$ and $b_0(\Gamma) = \dim H_0(\Gamma)$.  

A similar formula holds for modular graphs with legs. A leg can be
considered as an edge ending in a univalent vertex of a distinguished
kind, so, fix a linear operator $\zeta : V \to \setC$ and extend graphical
calculus so to evaluate univalent vertices to $\zeta$:
\begin{equation}\label{eq:modulargraphs}
    \int_V \exp \biggl\{\frac{1}{\hbar}
\left(\zeta(v)+\sum_{k=1}^\infty\sum_{g=0}^\infty
\hbar^g
\frac{S_{k,g}(v\tp{k})}{k!}
    \right)\biggr\} \ud\mu_\hbar(v)
    = \sum_{g,n=0}^\infty\hbar^g\sum_{\Gamma \in \textrm{M}(g,n)}
    \hbar^{-b_0(\Gamma)} \frac{Z(\Gamma)} {\card{\Aut\Gamma}}.
\end{equation}
\end{xmp}

\begin{xmp}[Asymptotic expansion of the ``free energy'' functional]
  The logarithm of the partition function $Z(x_*)$ is called the
  \emph{free energy} of the cyclic (resp. symmetric) algebra, and is
  denoted by $F(x_*)$; it admits a graph expansion, too. However,
  expansion of the partition function is a sum over \emph{all} graphs,
  whereas expansion of the free energy involves only
  \emph{connected} ones.

  \begin{lemma}\label{lemma:feynman-of-Z}
    The free energy \(F(x_*):= \log Z(x_*)\) admits a
    Feynman-Reshetikhin-Turaev expansion in ribbon (resp. ordinary)
    graphs, given by:
    \begin{equation}
      \label{eq:feynman-of-Z}
      F(x_*)
      = \sum_{\Gamma\text{\rmfamily\upshape\ connected}}
      \frac{Z_{x_*}(\Gamma)}{\card{\Aut\Gamma}}.
    \end{equation} 
  \end{lemma}
  \begin{proof}
    Exponentiate
    \begin{equation*}
      \sum_{\Gamma\textrm{ connected}}
      \frac{Z_{x_*}(\Gamma)}{\card{\Aut\Gamma}}
    \end{equation*}
    to find:
    \begin{equation*}
      \exp \left\{\sum_{\Gamma}
        \frac {Z_{x_*} (\Gamma)} {\card{\Aut\Gamma}}\right\}
      = \sum_{k=0}^\infty \oneof{k!} \sum_{\Gamma_1, \dots,
        \Gamma_k} \frac {Z_{x_*}(\Gamma_1) \cdots
        Z_{x_*}(\Gamma_k)} {\card{\Aut\Gamma_1} \cdots
        \card{\Aut\Gamma_k}},
    \end{equation*}
    where each $\Gamma_i$ is a connected graph.  Now recall that
    juxtaposition defines a tensor product $\otimes$ in the category of
    graphs (cf. \prettyref{dfn:graph-category}) and that $Z_{x_*}$ is
    multiplicative with respect to this structure:
    \begin{equation*}
      Z_{x_*} (\Gamma_1 \otimes \dots \otimes \Gamma_k)
      = Z_{x_*} (\Gamma_1) \otimes \dots \otimes Z_{x_*}
      (\Gamma_k).
    \end{equation*}
    Therefore,
    \begin{equation*} 
      \exp\left\{ \sum_{\Gamma\textrm{ connected}}
        \frac {Z_{x_*} (\Gamma)} {\card{\Aut\Gamma}}\right\} 
      = \sum_{k=0}^\infty 
      \sum_{\substack{\Gamma_1, \dots, \Gamma_k \\ \text{connected}}}
      \frac {Z_{x_*}(\Gamma_1 \otimes \cdots
        \otimes \Gamma_k)} {k! \card{\Aut\Gamma_1} \cdots
        \card{\Aut\Gamma_k}}
    \end{equation*}
    For a graph $\Phi$ having $k$ connected components $\Gamma_1, \dots,
    \Gamma_k$.  Let $I_\Phi$ be the set of all possible juxtapositions of
    $\Gamma_1, \ldots, \Gamma_k$; all graphs in $I_\Phi$ are isomorphic to $\Phi$. The
    semi-direct product $K$ of $\Perm{k}$ and $\Aut\Gamma_1 \times \dots \times
    \Aut\Gamma_k$ acts transitively on $I$; the stabilizer of any element
    is isomorphic to $\Aut\Phi$. Therefore,
    \begin{equation*}
      \card{I_\Phi} = \frac {\card{K}} {\card{\Stab\Phi}}
      = \frac {k! \card{\Aut\Gamma_1} \cdots \card{\Aut\Gamma_k}}
      {\card{\Aut\Phi}}.
    \end{equation*}
    So we reckon:
    \begin{equation*}
      \exp \left\{\sum_{\Gamma\textrm{ connected}}
        \frac {Z_{x_*} (\Gamma)} {\card{\Aut\Gamma}}\right\} 
      = \sum_{\Phi} \frac{Z_{x_*}(\Phi)}
      {\card{\Aut\Phi}}=Z(x_*).
    \end{equation*}
  \end{proof}
\end{xmp}


\everyxy={0,<2em,0em>:,(0,0.5),} 

\section{The 't~Hooft-Kontsevich model}
\label{sec:matrix-models}

This last section is devoted to the Kontsevich matrix model for 2D quantum
gravity, first defined in \cite{kontsevich;intersection-theory;1992}; it
embodies the ``standard matrix model'' of 't~Hooft as a particular case. A
cyclic algebra structure is introduced on the vector space $\Hermitian{N}$
of $N\times N$ Hermitian matrices; results from the previous sections apply.

Let \(V\) be an \(N\)-dimensional Hilbert space. The space \(\End(V)\)
has a natural Hermitian inner product
\begin{equation*} 
  \inner{X}{Y}:=\tr(X^*Y),
\end{equation*}
which induces the standard Euclidean inner product \(\inner{X}{Y} =
\tr(XY)\) on the real subspace
\begin{equation*}
  \Hermitian{N}:=\{X\in \End(V) | X^*=X\}
\end{equation*}
of Hermitian operators.  

For any positive definite Hermitian operator \(\Lambda\), we can define
a new Euclidean inner product on \(\Hermitian{N}\) by
\begin{equation*}
  \inner[\Lambda]{X}{Y} := 
  \onehalf \left(\tr(X\Lambda Y) + \tr(Y\Lambda X)\right).
\end{equation*}

Now, define cyclic tensors
\begin{equation*}
  T_k:\Hermitian{N}\tp{k} \ni X_1 \otimes X_2 \otimes \dots \otimes X_k 
  \mapsto 
  \tr(X_1X_2\cdots X_k) \in \setC;
\end{equation*}
These \(T_k\), together with the inner product $\inner[\Lambda]{-}{-}$,
define a cyclic algebra structure on \(\Hermitian{N}\) called the
Kontsevich model. A graphical calculus \(Z\) for the Kontsevich model is
defined on the PROP \(\RG\) of ribbon graphs.

\begin{lemma}\label{lem:Z1}
  The following formula holds:
  \begin{equation*}
    \label{eq:KM}
    Z (\Gamma) =
    \sum_{c} \prod_{\ell\in
      \Edges{\Gamma}} \frac{2}{\Lambda_{c(\ell^+)} + \Lambda_{c(\ell^-)}},
    \qquad c:\Holes{\Gamma} \to \{1, \dots, N\}
  \end{equation*}
  where \(\Lambda_1, \dots, \Lambda_N\) are the eigenvalues of \(\Lambda\), $c$ runs
  over all colorings of holes of $\Gamma$ in $N$ colors, and $\ell^+$,
  $\ell^-$ are the two holes bounded by the edges $\ell$ (they are not
  necessarily distinct).
\end{lemma}
\begin{proof}
  To evaluate \(Z(\Gamma)\) we need an explicit expression for the Casimir
  element \(\coev_{\Hermitian{N},\Lambda}(1)\) of the cyclic algebra \((
\Hermitian{N},
  \inner[\Lambda]{-}{-}, T_1, T_2, \dots)\).  Since \(\Lambda\) is Hermitian
  positive definite, there exists an orthonormal basis \(\{e_i\}\) of
  \(V\) in which
  \begin{equation*}
    \Lambda=\diag(\Lambda_1,\Lambda_2\dots,\Lambda_N),
  \end{equation*}
  for some \(\Lambda_i\) positive real numbers. Any choice of a like basis
  induces an identification of \(V\) with \(\setC^N\), and, consequently,
  of \(\End(V)\) with the space \(M_N(\setC)\) of \(N \times N\) complex
  matrices. Let \(\{E_{ij}\}\) be the canonical basis for \(M_N(\setC)\):
  \begin{equation*}
    (E_{ij})_{kl}=\delta_{ik}\delta_{jl}.
  \end{equation*}
  A basis for \(\Hermitian{N}\) is given by matrices
  \begin{equation*}
    e_{ij}=\begin{cases}
      \frac{1}{\sqrt{2}} (E_{ij}+E_{ji}) &\text{if \(i<j\),}\\
      E_{ii}                             &\text{if \(i=j\),}\\
      \sqrt{-\onehalf} (E_{ij}-E_{ji})   &\text{if \(i>j\).}
    \end{cases}
  \end{equation*}
  It is orthonormal with respect to the inner product
  \(\inner{-}{-}\), whereas 
  \begin{equation*}
    \inner[\Lambda]{e_{ij}}{e_{kl}}
    = \frac{\Lambda_i + \Lambda_j}{2} \delta_{ij,kl},
  \end{equation*}
  i.e., the matrix of \(\inner[\Lambda]{-}{-}\) with respect to the
  basis \(\{e_{ij}\}\) is
  \begin{equation*}
    g = \diag\left(\left\{ \frac{\Lambda_i + \Lambda_j}{2}
      \right\}\right).
  \end{equation*}
  So we get the following expression for the Casimir element:
  \begin{equation*}
    \coev_{\Hermitian{N},\Lambda}(1) = \sum_{i,j}
\frac{2}{\Lambda_i+\Lambda_j} e_{ij}
\otimes
    e_{ij}.
  \end{equation*}
  Rewrite this identity as:
  \begin{align*}
\coev_{\Hermitian{N},\Lambda}(1)&=\sum_{i}\frac{1}{{\Lambda_i}}e_{ii}\otimes
    e_{ii}+\sum_{i<j}\frac{2}{{\Lambda_i+\Lambda_j}}e_{ij}\otimes
    e_{ij}+\sum_{i>j}\frac{2}{{\Lambda_i+\Lambda_j}}e_{ij}\otimes
    e_{ij}\\
    &=\sum_{i}\frac{1}{{\Lambda_i}}e_{ii}\otimes
    e_{ii}+\sum_{i<j}\frac{2}{{\Lambda_i+\Lambda_j}}(e_{ij}\otimes
    e_{ij}+e_{ji}\otimes
    e_{ji}),
  \end{align*}
  but, for \(i<j\),
  \begin{align*}
    e_{ij}\otimes e_{ij} &+ e_{ji}\otimes e_{ji} = \\
    &\qquad + \frac{1}{
      2}(E_{ij}\otimes E_{ij} + E_{ij}\otimes E_{ji} + E_{ji}\otimes
    E_{ij}+E_{ji}\otimes E_{ji})\\
    &\qquad - \frac{1}{ 2}(E_{ij}\otimes E_{ij} - E_{ij}\otimes E_{ji} -
    E_{ji}\otimes
    E_{ij}+E_{ji}\otimes E_{ji})\\
    &= E_{ij}\otimes E_{ji} + E_{ji}\otimes E_{ij}.
  \end{align*}
  So,
  \begin{multline}\label{eq:casimir}
    \coev_{\Hermitian{N},\Lambda}(1)=
\sum_{i}\frac{1}{{\Lambda_i}}E_{ii}\otimes
    E_{ii}+\sum_{i<j}\frac{2}{{\Lambda_i+\Lambda_j}}(E_{ij}\otimes
    E_{ji}+E_{ji}\otimes
    E_{ij})\\
    = \sum_{i,j}\frac{2}{{\Lambda_i+\Lambda_j}}E_{ij}\otimes E_{ji}.
  \end{multline}

  According to  standard graphical calculus rules, evaluation
  \(Z(\Gamma)\) is performed through the correspondence
  \begin{equation*}
    {\xy*!LC\xybox{\rgvertex7\loose1\loose2\missing3%
        \loose4\loose5\missing6\loose7}\endxy}
    \leftrightarrow
    T_k,
    \qquad
    {\xy\vloop-\endxy}
    \leftrightarrow
    \coev_{\Hermitian{N},\Lambda}(1).
  \end{equation*}
  If we introduce the notation
  \begin{equation*}
    {\xy\vloop-,(0.05,0.5)*\txt{\({}_i\
        {}_j\)},(1,0.5)*\txt{\({}_j\
        {}_i\)}\endxy}=\frac{2}{\Lambda_i+\Lambda_j} E_{ij}\otimes E_{ji},
  \end{equation*}
  then we can depict \prettyref{eq:casimir} as
  \begin{equation*}
    {\xy\vloop-\endxy}
    = \sum_{i,j}
    {\xy\vloop-,(0.05,0.5)*\txt{\({}_i\
        {}_j\)},(1,0.5)*\txt{\({}_j\ {}_i\)}\endxy},
  \end{equation*}
  which turns \(Z(\Gamma)\) into a sum of ribbon graphs equipped with a
  number in \(\{1, \dots, N\}\) on each side of every edge, and
  operators \(T_k\) on each \(k\)-valent vertex.

  The map \(T_k\) is the restriction of a map \(T_k\) defined on
  \(M_N(\setC)\), namely, the trace of a \(k\)-fold product. We have
  \begin{equation}\label{eq:vertices}
    T_k(E_{i_1j_1}\otimes E_{i_2j_2}\otimes\cdots\otimes
    E_{i_{k}j_k})=\delta_{j_1i_2}\delta_{j_2i_3}\cdots
    \delta_{j_{k-1}i_k}\delta_{j_ki_1}.
  \end{equation}
  Therefore, the only graphs that give non-zero contribution to the sum
  giving $Z(\Gamma)$ are the ones whose boundary components have the same
  index on all edges --- that is, we need only account for graphs
  equipped with a map $c:\Holes{\Gamma} \to \{1, \dots, N\}$.  An edge whose
  sides are indexed $i,j$ brings in a factor $2/(\Lambda_i+\Lambda_j)$; combining
  this with \prettyref{eq:vertices}, we can conclude the proof.
\end{proof}

\begin{xmp}[The standard matrix model]
  Take $\Lambda = I$; formula \prettyref{eq:KM} specializes to 
  \begin{equation*}
    Z(\Gamma) = \sum_{c} \prod_{\ell\in \Edges{\Gamma}} \frac{2}{\Lambda_{c(\ell^+)} +
      \Lambda_{c(\ell^-)}}
    = \sum_c 1 = N^{\card{\Holes{\Gamma}}}.
  \end{equation*}
  Therefore, according to \prettyref{thm:FRT},
  \begin{equation}
    \int_{\Hermitian{N}} \exp \left\{ \frac{1}{\hbar}\sum_{j=1}^{\infty}
\frac{\tr X^j}{j}
    \right\} \ud\mu_I(X) = 
\sum_\Gamma \frac{N^{\card{\Holes{\Gamma}}}} {\card{\Aut
        \Gamma}}\hbar^{-\card{\Vertices{\Gamma}}}.
  \end{equation}
  This is known as the ``standard matrix model'' in physics
  literature.
\end{xmp}

\begin{xmp}[The orbifold Euler characteristics of $\M_g^n$]
  The Feynman-Reshetikhin-Turaev expansion formula for the
  Kontsevich model gives
  \begin{equation}\label{eq:FThighpot}
    \int_{\Hermitian{N}}
    \exp\left\{\sum_{j=1}^{\infty}x_j
      \frac{\tr X^{j}}{j}\right\}
    \ud\mu_\Lambda(X)=\sum_{\Gamma \in \R(0,0)}
    \frac{Z_{x_*}(\Gamma)}{\card{\Aut\Gamma}},
  \end{equation}
  where \(\R(0,0)\) is the set of isomorphism classes of (possibly
  disconnected) ribbon graphs of type (0,0). 

Now set  
  \begin{equation*}
    x_1=x_2=0, \quad 
    x_j=(\U)^j t^{\frac{j-2}{2}}, \quad
    \Lambda = I,
  \end{equation*}
  to find the graph enumeration formula \cite{bessis-itzykson-zuber;graphical-enumeration}:
  \begin{equation}
    \label{eq:enumeration}
    \log\int_{\Hermitian{N}}
    \exp\left\{\frac{1}{t}\sum_{j=3}^{\infty}(\sqrt{-t})^j
      \frac{\tr X^{j}}{j}\right\} \ud\mu_I(X) 
    = 
    \sum_{g,n} \sum_{\Gamma\in {\R}_g^n}
    \frac{(-1)^{\card{\Edges{\Gamma}}}} {\card{\Aut{\Gamma}}} 
    t^{2g-2+n} N^n,
  \end{equation}
  where ${\R}_g^n$ denotes the set of isomorphism classes of
  \emph{connected} genus $g$ ribbon graphs with $n$ boundary
  components.  It is well known (cf.
  \cite{harer;cohomology-of-moduli, mulase-penkava}) that the moduli
  space $\M_g^n$ of genus $g$ smooth complex curves with $n$
  (non-ordered) punctures has an orbifold triangulation whose cells
  are indexed by elements of $\R_g^n$ 
  Moreover, the local isotropy
  group of the cell $\Delta_\Gamma$ defined by $\Gamma$ is isomorphic to $\Aut
  \Gamma$, so the orbifold Euler characteristic of $\M_g^n$ can be
  computed to be
  \begin{equation*}
    \chi^{\textrm{orb}}(\M_g^n) = \sum_{\Gamma\in  
      {\R}_g^n}\frac{(-1)^{\dim \Delta_\Gamma}}{\card{\Aut{\Gamma}}},
  \end{equation*}
  see \cite{harer-zagier;euler-characteristic}. Since $\dim
  \Delta_\Gamma=\card{\Edges{\Gamma}}$, formula \prettyref{eq:enumeration} above
  is actually a generating series for the orbifold Euler
  characteristic of moduli spaces; indeed,
  \begin{equation*}
    \log \int_{\Hermitian{N}} \exp\left\{\frac{1}{t}\sum_{j=3}^\infty(\sqrt{-t})^j
      \frac{\tr 
        X^j}{j}\right\}\ud\mu(X)=\sum_{g,n}\chi^{\textrm{orb}}(\M_g^n)
    t^{2g-2+n}N^n.
  \end{equation*}
\end{xmp}


\appendix
\everyxy={/r24pt/:} 

\section{PROPs and Operads}
\label{appendix:PROPs}

The concept and structure of a PROP provides adequate language to
state properties of the graph families $\RT$, $\RG$, $\CG$ from a
categorical point of view. Since the definition of a PROP seems not to
be widely known, we recall here its axioms; our version is actually a
generalization of the one in \cite[pp.\ 
37--44]{adams;infinite-loop-spaces}, namely, we allow for the
category of indices to be an arbitrary monoidal one. We refer the
interested reader to \cite{adams;infinite-loop-spaces,
  boardman-vogt;homotopy-invariance, maclane;categorial-algebra,
  kelly;enriched-category} for a more thorough discussion and
background.

\begin{dfn}\label{dfn:PROP}
  A (non-symmetric) \emph{PROP} is the data of $(\catA_{\textrm{Ob}},
  ,\catA_{\textrm{Hom}},{\Pp}, \circ, \otimes, E,\break \Delta, j,
  \varphi^l, \varphi^r)$ where
  \begin{enumerate}[i)]  
  \item $\catA_{\textrm{Ob}}$ and $\catA_{\textrm{Hom}}$ are tensor
    categories, called respectively the \emph{category of objects} and
    the \emph{category of morphisms}; moreover $\catA_{\textrm{Hom}}$
    is symmetric;
  \item ${\Pp}:\left(\catA_{\textrm{Ob}}\right)^{\textrm{op}}\times 
    \catA_{\textrm{Ob}}\to\catA_{\textrm{Hom}}$ is a functor, called the
    \emph{\(\Hom\)-space} functor;
  \item $\circ_{A,B,C}:{\Pp}(B,C)\otimes {\Pp}(A,B)\to
    {\Pp}(A,C)$ is a natural
    transform, called the \emph{composition map};
  \item $\otimes_{A,B,C,D}:{\Pp}(A,B)\otimes {\Pp}(C,D)\to {\Pp}(A\otimes
    C,B\otimes D)$ is a natural
    transform;
  \item $(E,\Delta)$ is a co-associative co-algebra in
    $\catA_{\textrm{Hom}}$; 
  \item $j$ is a functorial map $j_A:E\to{\Pp}(A,A)$, 
    called the \emph{identity element};
  \item $\varphi^l_X:X\to E\otimes X$ and $\varphi^r_X:X\to X\otimes E$ are
    natural transformations; 
  \end{enumerate}
  such that all the diagrams of \prettyref{fig:PROP-axioms-1} and
  \prettyref{fig:PROP-axioms-2} commute.
\end{dfn}

We say that $(\catA_{\textrm{Ob}}, ,\catA_{\textrm{Hom}}, {\Pp},
\circ, \otimes, E, \Delta, j, \varphi^l, \varphi^r)$ is a
\(\catA_{\textrm{Hom}}\)-PROP over \(\catA_{\textrm{Ob}}\); we denote
it just by the symbol \(\Pp\).

\begin{rem}
  Unless the contrary is explicitly stated, \(E\) is the unit object
  of \(\catA_{\textrm{Hom}}\) and \(\Delta,\varphi^l,\varphi^r\) are defined in
  terms of the natural transformations \(l\) and \(r\) that are part
  of the tensor category structure on \(\catA_{\textrm{Hom}}\).  By
  abuse of language, elements of $\catA_{\textrm{Hom}}(E,{\Pp}(A,B))$
  are called \emph{elements} of ${\Pp}(A,B)$ and we write
  \(f\in{\Pp}(A,B)\) to mean \(f:E\to {\Pp}(A,B)\).
\end{rem}

The commutativity of \prettyref{eq:O1} means that composition of
morphisms in the PROP is associative.  Diagram \prettyref{eq:O1'}
expresses the condition
\begin{equation*}
  (f\otimes g) \circ (h\otimes k) = (f\circ h) \otimes (g\circ k).
\end{equation*}
Similarly, diagrams \prettyref{eq:OE} and \prettyref{eq:Oj} express
both the fact that $j_A$ acts as an identity and the compatibility
condition
\begin{equation*}
  j_A \otimes j_B = j_{A\otimes B}.
\end{equation*}
Finally, \prettyref{eq:O5} and~\prettyref{eq:O6} state that $E$ is an
identity element for the tensor product $\otimes$. 

The above requirements, admittedly cumbersome and obscure, may be
clarified by the following construction.  The coassociativity of the
comultiplication \(\Delta\) allows one to define a new category
$\catA_{\Pp}$ with:
\begin{enumerate}
\item the same objects as $\catA_{\textrm{Ob}}$;
\item morphisms given by
  \begin{equation*}
    \catA_{\Pp}(A,B):=\catA_{\textrm{Hom}}(E,{\Pp}(A,B)).
  \end{equation*}
\end{enumerate}
Composition of two morphisms \(f\in \catA_{\Pp}(B,C)\) and
\(g\in\catA_{\Pp}(A,B)\) is defined by
\begin{equation*}
  E\xrightarrow{\Delta} E\otimes E\xrightarrow{f\otimes
    g}{\Pp}(B,C)\otimes{\Pp}(A,B)\xrightarrow{\circ}{\Pp}(A,C).
\end{equation*}
The category \(\catA_{\Pp}\) is called the category \emph{underlying}
the PROP.  By the commutativity of \prettyref{eq:OE},
\(j_A\in{\Pp}(A,A)\) is the identity element of \(A \in \catA_{\Pp}\).
It is an easy exercise to verify that $\catA_{\Pp}$ is a monoidal
category with the tensor product $\otimes$. Note that every tensor
category \(\catA\) defines a \(\catVect\)-PROP by setting
\({\Pp}_{\catA}(X,Y):=\Hom(X,Y)\), for any two objects \(X\), \(Y\) of
\(\catA\). Since the unit object of \(\catVect\) is the base field
\(\fk\), the canonical isomorphism 
\begin{equation*}
\Hom(\fk,{\Pp}_{\catA}(X,Y))=\Hom(\fk,\Hom(X,Y))\simeq \Hom(X,Y)
\end{equation*}
gives a natural identification of \(\catA\) with the underlying category
of \({\Pp}_{\catA}\).

Maps \(\varphi^l\) and \(\varphi^r\) allow one to look at elements of
\({\Pp}(B,C)\) as left operators on \({\Pp}(A,B)\) (resp.  elements of
\({\Pp}(A,B)\) as right operators on \({\Pp}(B,C)\)) with values in
\({\Pp}(A,C)\). Indeed, the maps
\begin{equation*}
  \catA_{\textrm{Hom}}(E,{\Pp}(B,C))\to
  \catA_{\textrm{Hom}}({\Pp}(A,B),{\Pp}(A,C))
\end{equation*}
and
\begin{equation*}
  \catA_{\textrm{Hom}}(E,{\Pp}(A,B))\to
  \catA_{\textrm{Hom}}({\Pp}(B,C),{\Pp}(A,C))
\end{equation*}
are defined by
\begin{equation*} 
  {\Pp}(B,C)\ni f\mapsto\{{\Pp}(A,B)\xrightarrow{\varphi^l}
  E\otimes{\Pp}(A,B)\xrightarrow{f\otimes I}   
  {\Pp}(B,C)\otimes{\Pp}(A,B)\xrightarrow{\circ}   
  {\Pp}(A,C)\}
\end{equation*}
and
\begin{equation*}
  {\Pp}(A,B)\ni f\mapsto\{{\Pp}(B,C)\xrightarrow{\varphi^r}
  {\Pp}(B,C)\otimes E\xrightarrow{I\otimes f}
  {\Pp}(B,C)\otimes{\Pp}(A,B)\xrightarrow{\circ}
  {\Pp}(A,C)\}.
\end{equation*}
The commutativity of diagrams \prettyref{eq:O5} expresses
compatibility of this action with the composition of morphisms in the
PROP; commutativity of \prettyref{eq:O6} states that the identity
elements act trivially.

\begin{dfn} A braided (resp. symmetric) PROP is a PROP with the additional 
  datum of a natural morphism $\sigma_{A,B}:E\to{\Pp}(A\otimes B,B\otimes A)$
  inducing a braided (resp. symmetric) category structure on the
  category underlying the PROP.
\end{dfn}
\begin{xmp}
  Any braided (symmetric) tensor category is a braided
(symmetric) \(\catVect\)-PROP.  The category of Hilbert spaces is a
symmetric Banach spaces PROP.
  
  Moduli spaces of stable curves are an example of a symmetric
  algebraic-stacks-PROP over the monoidal category \(\catN\) of
  natural numbers.
\end{xmp}

\begin{dfn}
  A PROP $\Pp$ is \emph{linear} iff its category of \(\Hom\)-spaces is a
subcategory of \(\catVect\).
\end{dfn}
Let $\Pp$ be a linear PROP. 
\begin{dfn}
  Let $G$ be a subset of the set $\bigcup_{A, B} \Pp(A, B)$ of PROP
  operations.  The sub-PROP ${\Pp}_G$ \emph{generated} by $G$ is the
  smallest sub-PROP of $\Pp$ containing $G$.

  The linear PROP $\Pp$ is \emph{generated} by $G$ iff ${\Pp}_G = \Pp$.
\end{dfn}
 
\begin{dfn}
  Two PROPs ${\Pp}' := (\catA'_{\text{Ob}}, \catA'_{\text{Hom}}, \ldots)$
  and ${\Pp}'' := (\catA''_{\text{Ob}}, \catA''_{\text{Hom}}, \ldots)$ are
  deemed \emph{comparable} iff there exists a monoidal functor
  $\catA'_{\text{Hom}} \to \catA''_{\text{Hom}}$.
\end{dfn}
\begin{dfn}
  Given two comparable PROPs ${\Pp}'$ and ${\Pp}''$, define a morphism
  $\rho: \Pp' \to \Pp''$ to be a triple $(h, \rho_{\text{Ob}},
  \rho_{\text{Hom}})$, where:
  \begin{itemize}
  \item $h: \catA'_{\text{Hom}} \to \catA''_{\text{Hom}}$ is a monoidal
    functor,
  \item $\rho_{\text{Ob}}: \catA'_{\text{Ob}} \to \catA''_{\text{Ob}}$ is a
    tensor functor,
  \item $\rho_{\text{Hom}}$ is a natural transformation between the
    functors $h \circ {\Pp}'(-, -)$ and ${\Pp}''(\rho_{\text{Ob}}-,
    \rho_{\text{Ob}}-)$.
  \end{itemize}
  The triple $(h, \rho_{\text{Ob}}, \rho_{\text{Hom}})$ must satisfy
  conditions that express compatibility with the tensor structure on
  $\catA_{\text{Hom}}$.
   
  If $\rho: \Pp' \to \Pp''$ is an epimorphism, then we say that
  $\Pp''$ is a PROP quotient of $\Pp'$.
\end{dfn}

If $\catA'_{\textrm{Ob}} = \catA''_{\textrm{Ob}}$, $\catA'_{\textrm{Hom}}
=
\catA''_{\textrm{Hom}}$, $h=\Id$ and
$\rho_{\textrm{Ob}} = \Id$ , then PROP quotients are characterized by
kernels of maps $\rho_{\textrm{Hom}}(A,B): \Pp'(A,B) \to \Pp''(A,B)$.

Fix a vector space \(V\). 
\begin{dfn}
  The functor $\EndOp[V]: \catN\opp \times \catN \to \catVect$ defined by
  \begin{equation*}
    \EndOp[V](m,n) := \Hom(V^{\otimes m},V^{\otimes n}).
  \end{equation*}
  gives a structure of a \(\catVect\)-PROP over \(\catN\) to the
  monoidal category of tensor powers of $V$; it
  is called the \emph{endomorphism PROP} of \(V\).
\end{dfn}

\begin{dfn}
  Let $\Pp$ be a linear PROP over $\catN$. An action of $\Pp$ on a linear
  space $V$ is a morphism of PROPs $\Pp \to \EndOp[V]$, that is a
  collection of maps
  \begin{equation*}
    \Pp(p,q) \to \Hom(V\tp{p}, V\tp{q}),
  \end{equation*}
  satisfying obvious compatibility conditions. A representation of \(\Pp\)
is a pair \((V,\rho)\) where \(V\) is a vector space and \(\rho\) is an
action of \(\Pp\) on \(V\). If  \((V,\rho)\) is a representation of
\(\Pp\), then \(V\) is called a \(\Pp\)-algebra. 
\end{dfn}

\begin{rem} PROPs are deeply related to operads,
see \cite{adams;infinite-loop-spaces, may;loop-spaces}.
  Indeed, if \(\Pp\) is a symmetric \(\catA_{\textrm{Hom}}\)-PROP
  over \(\catN\), the collection
\begin{equation*}
  {\Oo}_{\Pp}(n):={\Pp}(n,1)
\end{equation*}
is a \(\catA_{\textrm{Hom}}\)-operad. Conversely, if
the collection \({\Oo}(n)\) is a May operad,
\begin{equation*}
{\Pp}_{\Oo}(m,l):=\bigoplus_{\substack{m_1,m_2,\dots,m_l\\
\sum_i m_i=m}}
{\Pp}(m_1)\otimes\cdots\otimes{\Pp}(m_l)
\end{equation*}
defines a symmetric PROP. Notice that
\({\Oo}_{{\Pp}_{\Oo}}={\Oo}\), whereas it is only
\({\Pp}_{{\Oo}_{\Pp}}\subseteq{\Pp}\)
in general. 
\end{rem}

\begin{sidewaystable}[p]
  \label{fig:PROP-axioms-1}
  \caption{Diagrams expressing compatibility relations of structure
    maps in a PROP, part I.}
  \begin{align}
    \label{eq:O1}\tag{O1}
    &{\xymatrix@C=-6pt{
        \displaystyle{{\Pp}(C,D)\otimes({\Pp}(B,C)\otimes{\Pp}(A,B))}
        & & \displaystyle{({\Pp}(C,D)\otimes {\Pp}(B,C))\otimes
          {\Pp}(A,B)}
        \\
        \displaystyle{{\Pp}(C,D)\otimes {\Pp}(A,C)} & &
        \displaystyle{{\Pp}(B,D)\otimes {\Pp}(A,B)}
        \\
        &\displaystyle{{\Pp}(A,D)}& \ar "1,3";"1,1"_{a} \ar
        "1,1";"2,1"_{\Id\otimes\circ} \ar
        "1,3";"2,3"^{\circ\otimes\Id} \ar "2,1";"3,2"_{\circ} \ar
        "2,3";"3,2"^{\circ} }}
  \end{align}
  \begin{align}
    \label{eq:O1'}\tag{O2}
    {\xymatrix@C=-24pt{ \displaystyle{({\Pp}(A,B)\otimes {\Pp}(C,D))
          \otimes ({\Pp}(F,A) \otimes \Pp(G,C))} & &
        \displaystyle{({\Pp}(A,B)\otimes {\Pp}(F,A)) \otimes
          ({\Pp}(C,D) \otimes \Pp(G,C))}
        \\
        \displaystyle{{\Pp}(A \otimes C, B \otimes D) \otimes {\Pp}(F
          \otimes G, A \otimes C)} & & \displaystyle{{\Pp}(F,B)\otimes
          {\Pp}(G,D)}
        \\
        &\displaystyle{{\Pp}(F\otimes G, B\otimes D)}& \ar@{<->}
        "1,3";"1,1"_{\textrm{shuffle}} \ar "1,1";"2,1"_{(-\otimes-)
          \otimes (-\otimes-)} \ar "1,3";"2,3"^{\circ\otimes\circ} \ar
        "2,1";"3,2"_{\circ} \ar "2,3";"3,2"^{\otimes} }}
  \end{align}
\end{sidewaystable}
\begin{table}[p]
  \label{fig:PROP-axioms-2}
  \caption{Diagrams expressing compatibility relations of structure
    maps in a PROP, part II.}
  \begin{align}
    \label{eq:OE}\tag{O3}
    {\xymatrix@C=18pt{
        \displaystyle{E}& \displaystyle{E\otimes E}\\
        \displaystyle{\Pp(A,B)} & \displaystyle{\Pp(A,B)\otimes
          \Pp(A,A)} \ar "1,1";"1,2"^{\Delta} \ar "1,1";"2,1"_{\forall
          f} \ar "1,2";"2,2"^{f \otimes j_A} \ar
        "2,2";"2,1"_(.60){\circ} }} 
    \quad
    &{\xymatrix@C=18pt{
        \displaystyle{E}& \displaystyle{E\otimes E}\\
        \displaystyle{\Pp(A,B)} & \displaystyle{\Pp(B,B)\otimes
          \Pp(A,B)} \ar "1,1";"1,2"^{\Delta} \ar "1,1";"2,1"_{\forall
          f} \ar "1,2";"2,2"^{j_B \otimes f} \ar
        "2,2";"2,1"_(.60){\circ} }}
    \\[12pt]
    \label{eq:Oj}\tag{O4}
    {\xymatrix@C=18pt{%
        \displaystyle{E} \otimes \displaystyle{E} & \displaystyle{E} \\
        \displaystyle{\Pp(A,A) \otimes \Pp(B,B)} & \displaystyle{\Pp
          (A\otimes B, A\otimes B)} \ar "1,2";"1,1"_{\varphi^r} \ar
        "1,1";"2,1"_{j_A \otimes j_B} \ar "1,2";"2,2"^{j_{A\otimes B}}
        \ar "2,1";"2,2"_{\otimes} }}
    \quad
    &{\xymatrix@C=18pt{%
        \displaystyle{E} \otimes \displaystyle{E} & \displaystyle{E} \\
        \displaystyle{\Pp(A,A) \otimes \Pp(B,B)} & \displaystyle{\Pp
          (A\otimes B, A\otimes B)} \ar "1,1";"1,2"^{\varphi^l} \ar
        "1,1";"2,1"_{j_A \otimes j_B} \ar "1,2";"2,2"^{j_{A\otimes B}}
        \ar "2,1";"2,2"_{\otimes} }}
    \\[12pt]
    \label{eq:O5}\tag{O5}
    {\xymatrix@C=48pt{
        \displaystyle{X}& \displaystyle{E\otimes X}\\
        \displaystyle{E\otimes X} & \displaystyle{E\otimes E\otimes X}
        \ar "1,1";"1,2"^{\varphi^l_X} \ar "1,1";"2,1"_{\varphi^l_X}
        \ar "1,2";"2,2"^{\varphi^l_{E\otimes X}} \ar
        "2,1";"2,2"^{\Delta\otimes \Id_X} }} 
    \quad
    &{\xymatrix@C=48pt{
        \displaystyle{X}& \displaystyle{X\otimes E}\\
        \displaystyle{X\otimes E} & \displaystyle{X\otimes E\otimes E}
        \ar "1,1";"1,2"^{\varphi^r_X} \ar "1,1";"2,1"_{\varphi^r_X}
        \ar "1,2";"2,2"^{\varphi^r_{X\otimes E}} \ar
        "2,1";"2,2"^{\Id_X\otimes \Delta} }}
    \\[12pt]
    \label{eq:O6}\tag{O6}
    {\xymatrix@C=12pt{
        \displaystyle{\Pp(A,B)}& \displaystyle{E\otimes \Pp(A,B)}\\
        \displaystyle{\Pp(A,B)} & \displaystyle{\Pp(B,B)\otimes
          \Pp(A,B)} \ar "1,1";"1,2"^{\varphi^l} \ar "1,1";"2,1"_{\Id}
        \ar "1,2";"2,2"^{j_B\otimes \Id} \ar "2,2";"2,1"_(.60){\circ}
      }} 
    \quad
    &{\xymatrix@C=12pt{
        \displaystyle{\Pp(A,B)}& \displaystyle{\Pp(A,B)\otimes E}\\
        \displaystyle{\Pp(A,B)} & \displaystyle{\Pp(A,B)\otimes
          \Pp(A,A)} \ar "1,1";"1,2"^{\varphi^r} \ar "1,1";"2,1"_{\Id}
        \ar "1,2";"2,2"^{\Id\otimes j_A} \ar "2,2";"2,1"_(.60){\circ}
      }}
  \end{align}
\end{table}


\nonumsection{Acknowledgements}
We are pleased to thank our advisor E. Arbarello, for having
stimulated our interest in this subject; we are deeply indebted to G.
Bini and F. Gavarini, who daringly read a very early manuscript of
this paper, and made a number of helpful suggestions.  We wish to
thank R. Oeckl and J. Stasheff for very important remarks.

\nonumsection{References}
\bibliography{math.AG}
\bibliographystyle{hunsrt}

\end{document}